\documentclass[ssy, 12pt]{imsart}

\usepackage{amsfonts,amsthm,amsmath,amssymb,epsfig,epstopdf,epic}
\RequirePackage[colorlinks,citecolor=blue,urlcolor=blue]{hyperref}

\arxiv{arXiv:0000.0000}

\startlocaldefs  
\numberwithin{equation}{section}
\theoremstyle{plain}
\newtheorem{them}{Theorem}

\newtheorem{lema}{Lemma}
\newtheorem{prop}{Proposition}
\newtheorem{defn}{Definition}
\newtheorem{remk}{Remark}

\endlocaldefs

\usepackage{geometry}

\begin{document}

\begin{frontmatter}

\title{Stochastic Monotonicity of Markovian Multi-class Queueing Networks}
\runtitle{Stochastic Monotonicity of McQNs}


\begin{aug}
  \author{\fnms{H.} \snm{Leahu}\corref{}\ead[label=e1]{haralambie@gmail.com}},
  \author{\fnms{M.} \snm{Mandjes}\ead[label=e2]{m.r.h.mandjes@uva.nl}}

  \runauthor{H. Leahu and M. Mandjes}

  \affiliation{University of Amsterdam}

  \address{University of Amsterdam, Science Park 904, 1098 XH, Amsterdam, The Netherlands; \\ \printead{e1,e2}}
 
\end{aug}

\begin{abstract}
\noindent Multi-class queueing networks (McQNs) extend the classical concept of Jackson network by allowing jobs of different classes to visit the same server. While such a generalization seems rather natural, from a structural perspective there is a significant gap between the two concepts. Nice analytical features of Jackson networks, such as stability conditions, product-form equilibrium distributions, and stochastic monotonicity do not immediately carry over to the multi-class framework.

\noindent The aim of this paper is to shed some light on this structural gap, focusing on monotonicity properties. To this end, we introduce and study a class of Markov processes, which we call \emph{Q-processes}, modeling the time evolution of the network configuration of any open, work-conservative McQN having exponential service times and {Poisson input}. We define a new monotonicity notion tailored for this class of processes. Our main result is that we show monotonicity for a large class of McQN models, covering virtually all instances of practical interest. This leads to interesting properties which are commonly encountered for `traditional' queueing processes, such as (i)~monotonicity with respect to external arrival rates and (ii)~star-convexity of the stability region (with respect to the external arrival rates); such properties are well known for Jackson networks, but had not been established
at this level of generality. 

\noindent This research was partly motivated by the recent development of a simulation-based method which allows one to
numerically determine the stability region of a McQN parametrized in terms of the arrival rates vector.
\end{abstract}



\end{frontmatter}


\newpage

\section{Introduction}\label{intro} 

Multi-class queueing networks (McQNs) arise as natural generalizations of conventional Jackson networks: while in Jackson networks
each station (server) acts as a $\cdot/M/1$ {\it single-class} queue, in McQNs each network station is a {\it multi-class} queue.
McQNs are particularly suitable for describing complex manufacturing systems (to be thought of as assembly lines) as they allow jobs 
(or, in queuing lingo, customers) visiting multiple times the same station to have different service requirements and/or a different routing scheme. One can think of situations in which a piece entering the system undergoes some physical transformation during the process, hence its processing time at a given station (as well as its next destination) might depend on the processing stage. Other possible applications include packet transmission models in telecommunication networks and distributed systems in computer science, where packets/tasks of different types can be routed to the same server in order to control resource utilization.
\vspace{2mm}

\noindent \emph{Motivation and Background}: The above generalization involves a number of complications. For instance, the following
mathematical challenges arise:
\begin{enumerate}
\item [(I)] Although the network-configuration process is expected to be Markovian (provided that all service and exogenous inter-arrival times are independent, exponentially distributed), the Markovian structure (the state-space and the transition dynamics) of a McQN is by no means straightforward, as it depends on the service disciplines of the underlying stations. While for some particular
service disciplines the structure becomes quite simple, a unified (Markovian) framework is lacking and this makes the analysis of such
networks rather difficult. The Markovian modeling of McQNs is important for both theoretical and practical reasons. On one hand it allows one to use the powerful Markov process machinery to derive analytical properties of the underlying network-configuration process, while on the other hand it facilitates the use of standard simulation methods for such networks.
\item [(II)] Stochastic monotonicity is, in general, a desirable property which is widely used in applications, e.g.\ optimization. Jackson networks satisfy the following monotonicity condition: if one increases any flow of jobs entering the network 
then the resulting queues (at each station, at any given time) will not decrease (in a stochastic sense) \cite{S-Yao}. However, a straightforward generalization of this result to the multi-class framework is impeded by the inhomogeneous nature of the  queues. 
More specifically, while homogeneous queues can be naturally ordered with respect to their length, it is not clear what ordering for inhomogeneous queues (consisting of different types of jobs) would be suitable for extending the stochastic monotonicity concept
to the multi-class setup.
\item [(III)] Stability is arguably a crucial (asymptotic) property of a queueing network. In this context, stability refers to positive recurrence of the associated Markov process. For Jackson networks, stability is equivalent to subcriticality, i.e.\ traffic rate below 1 at every queue. In the multi-class framework, however, this is not the case, as illustrated by various counter examples; see \cite{RS:93,Dai:95,Bramson:94,KS:90,S:94}. As such, stability conditions are not available in closed form, in general. 

\noindent  In addition, stability is closely related to monotonicity properties. More specifically, validity of monotonicity properties is expected to imply that stability of a queueing network is a monotone property. This is obviously the case for  Jackson networks, as stability (i.e. subcriticality) is a monotone property with respect to both (external) arrival, service and traffic rates. For McQNs, however, some results in the literature indicate that this is not necessarily the case, when parametrized with respect to service rates \cite{Dai:99,Bramson:95} or  {traffic rates} \cite{Dumas}. This naturally raises the question (but also casts some doubts on) whether stability of a McQN is monotone with respect to (external) arrival rates.
\end{enumerate}

(I) has been addressed in \cite{Bramson:08} where a Markovian formalism has been proposed in a fairly general framework (with no restrictions over the underlying distributions).   

(II) is a classical topic in applied probability and has attracted considerable attention over the past decades. We refer to
\cite{S-Yao} for monotonicity results pertaining to queueing networks and \cite{Massey:87} (the references therein) for a standard theory tailored to Markov processes. Neither of the approaches, however, covers the multi-class setup.  

Finally, much research has been invested in the 1990's into (III). The most successful approach to studying the stability of McQNs
is based on fluid (model) limits, an asymptotic technique originally introduced in \cite{RS:93} and further expanded in \cite{Dai:95}
which relates the stability of a McQN to that of its associated fluid model. Namely, stability of the fluid model implies that of the McQN; see \cite[Section 5.5]{Bramson:08} for more background, including counter-examples showing that the converse is not true. 

Stochastic monotonicity of (Markovian) McQNs became relevant recently, when, motivated by the lack of closed-form stability conditions, the authors investigated simulation-based methods for approximating stability regions \cite{LM:2016}. Extensive simulation experiments indicated that a certain form of (stochastic) monotonicity with respect to external arrival rates could still be expected, even for McQNs for which no stability conditions are known. Such a property would be enough to guarantee, for instance, that stability is a monotone property with respect to arrival rates.\vspace{2mm}


\noindent
\emph{Contributions}: In this paper, we restrict our analysis to {\it Markovian} McQNs. To this end, we introduce a new class of Markovian processes, called Q-processes, formalizing the network-configuration process associated with a Markovian McQN. Concerning questions (I) -- (III) above, the main contributions of this paper are the following: 
\begin{enumerate}
\item [1.] We develop a novel stochastic monotonicity concept tailored to Q-processes, called $\mathcal{F}$-monotonicity, and identify
a set of rather general conditions (fulfilled by virtually all McQNs of  practical interest) which guarantee $\mathcal{F}$-monotonicity; see Theorem \ref{HQ:them}. Furthermore, we show that $\mathcal{F}$-monotonicity implies monotonic behavior with respect to both (external) arrival rates and (when started empty) with time, extending in this way the well known results from Jackson networks.
\item [2.] Secondly, we prove that for $\mathcal{F}$-monotone Q-processes stability is a monotone property with respect to external arrival rates; see Theorem \ref{regio:them}. In particular, the stability region is an open, star-shaped domain, having the origin as
a vantage point. This result formally validates the numerical findings in \cite{LM:2016}.
\end{enumerate}

\noindent\emph{Approach}: Restricting our analysis to {\it Markovian} McQNs is necessary for several reasons. In the first place, it makes the corresponding network-configuration process Markovian, which significantly simplifies the formalism. Secondly, deviating from the Markovian framework leads to serious complications, as stability \cite{Dai:04} and monotonicity \cite{Whitt:93} are sensitive to changing the shape of the underlying distributions.  

A prerequisite for our investigation is an appropriate modeling framework that is general enough to cover a wide range of Markovian queue-related processes. More specifically, we need to define a state space which is wide enough to accommodate various types of
queue configurations together with an underlying structure which allows one to define the (Markovian) transitions performed by the network configuration process. We formalize the queue at a given station as an ordered sequence of digits, where the class of each job in the queue is represented by a specific digit; the order of the digits is interpreted as the order in which they are due to receive service, as long as no other arrivals occur in the queue. Furthermore, to accommodate various queueing disciplines, e.g.\ priority rules, we need to introduce a family of \emph{insertion operators} indicating how an extra-digit (representing the class of a new job arriving in the queue) will be placed in the sequence. Finally, in order to cover processor-sharing disciplines, we shall also introduce the concept of \emph{service allocation}, i.e.\ a probability distribution on the set of digits specifying the fraction of service allocated by the server to each class present in the queue. The space of ordered sequences endowed with a family of insertion operators and a service allocation will be called a \emph{space of multi-class configurations} and will be the building block for defining a Q-process; these facts are formalized in Section \ref{q:sec}.

Furthermore, we introduce the concept of $\mathcal{F}$\emph{-monotonicity}, which is a weaker form of stochastic monotonicity for Markov 
processes on ordered spaces. To be more specific, consider a Markov process on a (partially) ordered space. The ordering
induces a class of (real-valued) increasing functions and a stochastic ordering on the probability distributions on the underlying
state space \cite{Massey:87}. Standard stochastic monotonicity, cf.\ \cite{Massey:87}, presumes that two versions of the process
started in a pair of ordered states remain (stochastically) ordered at any time. Such a property does not hold, in general, for Q-processes. On the other hand, $\mathcal{F}$-monotonicity, for some given (sub)class of increasing functions $\mathcal{F}$, still presumes an ordering between the two versions, but with respect to a different (weaker) stochastic ordering, determined by the subclass $\mathcal{F}$. Although $\mathcal{F}$\emph{-monotonicity} is weaker than the standard stochastic monotonicity for Markov processes
on ordered spaces and requires essentially different proof techniques, it still retains most of the relevant properties of stochastic monotonicity and, as it turns out, it is better suited in the context of Q-processes.

Finally, we consider the vector of (external) arrival rates as a parameter of the Q-process corresponding to a given McQN and express the corresponding stability region (the set of parameters which make the process stable) as the support of some limiting functional of the process, defined on the parameter space. Based on analytical properties of this functional, we derive relevant properties of the corresponding stability region.\vspace{2mm}

\noindent \emph{Organization of the paper}: The paper is organized as follows. In Section \ref{mcqn:sec} we provide a brief account
of the mathematical concept of McQN. In Section \ref{q:sec} we define the concept of Q-process, i.e.\ the general stochastic process
model for the dynamics of a Markovian McQN. Then, in Section \ref{mon:sec} we introduce and elaborate on the (novel) concept of $\mathcal{F}$-monotonicity, while in Section \ref{stab:sec} we investigate stability properties of $\mathcal{F}$-monotone Q-processes. Finally, in Section \ref{sec:concl} we illustrate the practical importance of our results, by pointing out their relevance in developing numerical methods for evaluating the stability region associated with a Q-process (resp. McQN).\vspace{2mm}

\noindent \emph{Notations and Conventions}: In this paper, we shall employ the following notation. In the first place,  by $\mathbb{N}$ we
denote the set of non-negative integers $\mathbb{N}:=\{0,1,2,\ldots\}$ and by $\mathbb{R}$ the set of real
numbers. For a denumerable set $\mathcal{J}$ and $\mathbb{S}=\mathbb{N},\mathbb{R}$ we denote by $\mathbb{S}^\mathcal{J}$ the set
of $\mathcal{J}$-labeled vectors over $\mathbb{S}$; alternatively, $\mathbb{S}^\mathcal{J}$ defines the space of all mappings
$u:\mathcal{J}\longrightarrow\mathbb{S}$. When $\mathcal{J}=\{1,\ldots,d\}$ we  use the simplified notation $\mathbb{S}^d$.
Moreover, $\mathbf{0}$  denotes the null vector $(0,\ldots,0)$ and $\mathbf{1}_\Gamma$, for $\Gamma\subseteq\mathcal{J}$,  
denotes the characteristic vector (mapping) of $\Gamma$, defined via 
$$(\mathbf{1}_\Gamma)_\jmath=
\begin{cases}
1,\:\jmath\in\Gamma, \\
0,\:\jmath\notin\Gamma.
\end{cases}$$
When $\Gamma=\mathcal{J}$ we use the simplified notation $\mathbf{1}$ instead $\mathbf{1}_\Gamma$.

For an arbitrary vector ${\boldsymbol x}=(x_\jmath:\jmath\in\mathcal{J})\in\mathbb{S}^\mathcal{J}$ we define the $1$-norm 
$$\|{\boldsymbol x}\|:=\sum_{\jmath\in\mathcal{J}}|x_\jmath|,$$ 
whenever the r.h.s.\ is finite, and we extend this notation to bounded linear operators ${\boldsymbol U}$ defined on (subspaces of) $\mathbb{S}^\mathcal{J}$, viz. $\|{\boldsymbol U}\|:=\sup\{|{\boldsymbol U}{\boldsymbol x}|:\|{\boldsymbol x}\|\leq 1\}$. 
On $\mathbb{S}^\mathcal{J}$ we denote the natural (componentwise) ordering ${\boldsymbol x}\leq{\boldsymbol z}$, resp.\
${\boldsymbol x}<{\boldsymbol z}$, if $x_\jmath\leq z_\jmath$, resp.\ $x_\jmath < z_\jmath$, for $\jmath\in\mathcal{J}$.
Finally, $\varsigma[{\boldsymbol x}]:=\{\jmath\in\mathcal{J}:x_\jmath\neq 0\}$ denotes the support of
${\boldsymbol x}\in\mathbb{S}^\mathcal{J}$ and $\boldsymbol{I}$ denotes the identity operator on $\mathbb{S}^\mathcal{J}$.

 
\newpage 

\section{Markovian Multi-class Queueing Networks}\label{mcqn:sec}

Following the exposition in \cite{Dai:95}, we consider a multi-server network, comprising $\aleph\geq 1$ single servers, labeled $1,\ldots,\aleph$. The network is used by $d\geq 1$ classes of jobs in such a way that each class $k$ job, at any time, requires
service at a fixed server, denoted by $S(k)$. Once the service at $S(k)$ is finished, it either becomes a job of class $l$  with probability $R_{kl}$ (independently of all routing history), or leaves the system with (exit) probability $$R_{k0}:=1-\sum_{l=1}^d R_{kl}.$$ The routing matrix $\boldsymbol{R}=\{R_{kl}\}_{k,l=1,\ldots,d}$ is assumed transient (substochastic), i.e.,
\begin{equation}\label{R-series:eq}
\boldsymbol{I}+\boldsymbol{R}+\boldsymbol{R}^2+\ldots=(\boldsymbol{I}-\boldsymbol{R})^{-1}.
\end{equation}
This guarantees that any job entering the network visits a finite number of classes
before leaving the network, almost surely; in standard queueing language, the network is said to be \emph{open}.

Each class $k$ has its own exogenous (possibly null) arrival stream, regulated by a Poisson process with rate $\theta_k\geq0$ and requires i.i.d.\ service times, exponentially distributed with rate $\beta_k>0$, independent of everything else; a null arrival process corresponds to a class with no external input, which models, for instance, intermediate processing stages of a certain class. Note that any class is identified with its server, its specific routing probabilities, its specific exogenous arrival process and its specific service-time distribution.

We let $\mathcal{K}_i:=S^{-1}(i)$ denote the set of classes served by station $i$ and assume (without loss of generality)
that $\mathcal{K}_i\neq\emptyset$, for all $i$ (equivalently, the mapping $S$ is surjective), i.e., any server is used. In particular, it holds that $\aleph\leq d$.

Finally, each server employs its own non-idling service discipline (i.e., different servers may have different service disciplines)
and has infinite buffer capacity. 

The above introduced McQN concept extends many known classes of queueing networks. For instance, when the mapping $S$ is bijective (in particular, $\aleph=d$) one obtains a \emph{Jackson network}. Moreover, if the service rate for any class $k$ only
depends on the server $S(k)$ then we recover the concept of \emph{Kelly network};
see \cite{Kelly}. Finally, if there exists only one class with a non-null exogenous arrival process and all jobs have the same (deterministic) routing, visiting all classes exactly once (in the same order), then the network is called a \emph{re-entrant line}. Re-entrant lines provide  popular instances of McQNs, as they can be used to model  (assembly) manufacturing lines.
 
\section{Q-processes}\label{q:sec}

The objective of this section is to introduce the concept of Q-process, a Markovian process modeling the network configuration
dynamics (over time) of a Markovian McQN, as described in Section \ref{mcqn:sec}. In order to construct a model general enough to accommodate all the usual service disciplines, we need a rather intricate formalism which we briefly introduce below. Technical
details and definitions can be found in Appendix \ref{space:sec}.\vspace{2mm}

The building blocks of the Q-process model are the \emph{spaces of multi-class configurations} which formalize the basic structure
required for modeling the dynamics of a single-server, multi-class queue. More specifically, we consider a finite set $\mathcal{K}$
of job classes to be processed by a server. Since jobs of the same class are (probabilistically) exchangeable, the class of all
possible queue configurations is modeled by the augmented space $\overline{\mathbb{Q}}[\mathcal{K}]=\mathbb{Q}[\mathcal{K}]\cup\{\emptyset\}$, where
\begin{equation}\label{smcc:eq}
\mathbb{Q}[\mathcal{K}]:=\left\{p=(k_1,\ldots,k_n):\:n\geq 1,\:k_1,\ldots,k_n\in\mathcal{K}\right\},
\end{equation} 
denotes the space of all finite (ordered) sequences with elements in $\mathcal{K}$. 

Furthermore, in order to formalize the queue configuration dynamics, we need to introduce some additional structure on $\overline{\mathbb{Q}}[\mathcal{K}]$, modeling arrival and departure events. More specifically, we need:
\begin{itemize}
\item Insertion operators $\mathcal{I}_k:\mathbb{Q}[\mathcal{K}]\longrightarrow\mathbb{Q}[\mathcal{K}]$ which specify how an
incoming job of class $k\in\mathcal{K}$ will be placed in the (new) queue configuration, i.e., there exists a representation
(concatenation) $p=(p',p'')$ such that $\mathcal{I}_k(p)=(p',k,p'')$; both $p'$ and $p''$ are allowed to be $\emptyset$, but $p''$ may
not contain any $k$-digits (jobs in the same class may not overtake each other). The family of all insertion operators $\{\mathcal{I}_k:k\in\mathcal{K}\}$ defines a \emph{queue policy} on $\mathbb{Q}[\mathcal{K}]$.
\item Deletion operators $\mathcal{D}_k:\mathbb{Q}[\mathcal{K}]\longrightarrow\overline{\mathbb{Q}}[\mathcal{K}]$ which specify
how a job of class $k\in\mathcal{K}$ is removed from the queue configuration.
\item A service allocation (mapping) $\mathcal{V}:=(\mathcal{V}_k:k\in\mathcal{K}):\mathbb{Q}[\mathcal{K}]\longrightarrow\mathcal{P}[\mathcal{K}]$,
where $\mathcal{P}[\mathcal{K}]$ denotes the simplex of probability vectors $\boldsymbol{w}:=(\mathcal{V}_k:k\in\mathcal{K})$ on $\mathcal{K}$, specifying which class(es) receive service and the corresponding fraction of server capacity, in a given (nonempty) configuration. Since $\mathcal{V}_k(p)$ denotes the fraction of server capacity allocated to class $k$, we impose that $\mathcal{V}_k(p)=0$ if there is no $k$-digit
in the configuration $p$. 
\end{itemize}

The above elements relate to a McQN model as follows. The insertion operators are set in accordance with the queue discipline of the server. The order of the digits in a sequence $p$ does not necessarily reflect the order of arrivals, but rather the order in which jobs will be considered for service. Furthermore, a deletion operator
$\mathcal{D}_k$ removes (by convention) the first $k$-digit from the sequence; if no such digit exists, it leaves the sequence unchanged. Finally, a service allocation distributes the server capacity among the jobs/classes present in the queue; that is, each class receives a given fraction of server capacity (which, by convention,
is assigned to its first representative) in accordance with the service discipline of the server.

Note that in this modeling paradigm the service discipline of a particular server is identified with a combination of queue policy
and service allocation, which, in principle, can be chosen independently from each other. Note, however, that in some cases there are
multiple ways to model a particular service discipline; for instance, when the service allocation is order insensitive (which is typically the case for processor sharing disciplines), the queue policy becomes irrelevant (a default can be used).
For specific examples of queue policies and service allocations and how do they relate to the usual service disciplines, see Appendix \ref{space:sec}.

\begin{remk}
\emph{The mappings $\mathcal{I}_k$, $\mathcal{D}_k$ and $\mathcal{V}_k$ (which are defined on $\mathbb{Q}[\mathcal{K}]$) admit natural extensions to the augmented space $\overline{\mathbb{Q}}[\mathcal{K}]$, as follows: $\mathcal{I}_k(\emptyset):=(k)$, $\mathcal{D}_k(\emptyset):=\emptyset$ and $\mathcal{V}_k(\emptyset):=0$;
note, however, that $\mathcal{V}(\emptyset)=\mathbf{0}\notin\mathcal{P}[\mathcal{K}]$}.\hfill$\diamond$
\end{remk}

\begin{defn} \emph{A \emph{space of multi-class configurations} over $\mathcal{K}$ is the space
$$\left\{\overline{\mathbb{Q}}[\mathcal{K}];\mathcal{I}_k,\mathcal{D}_k,\mathcal{V}_k:k\in\mathcal{K}\right\},$$
of finite ordered sequences over $\mathcal{K}$, endowed with a family of insertion/deletion operators 
$\{\mathcal{I}_k/\mathcal{D}_k:k\in\mathcal{K}\}$ and service allocation $\mathcal{V}=(\mathcal{V}_k:k\in\mathcal{K})$}.\hfill$\diamond$
\end{defn}

Let $1\leq\aleph\leq d$ and assume that $\{\mathcal{K}_i\}_{i=1}^\aleph$ is a partition of the set $\{1,\ldots,d\}$.
Furthermore, we assume that 
$\left\{\overline{\mathbb{Q}}[\mathcal{K}_i];\mathcal{I}_k,\mathcal{D}_k,\mathcal{V}_k:k\in\mathcal{K}_i\right\}$ is a space of multi-class configurations over $\mathcal{K}_i$, for $i=1,\ldots,\aleph$ and let
\begin{equation}\label{state-sp:eq}
\mathbb{X}: = \overline{\mathbb{Q}}[\mathcal{K}_1]\times\ldots\times\overline{\mathbb{Q}}[\mathcal{K}_\aleph]
= \left\{[p_1,\ldots,p_\aleph],\:p_i\in\overline{\mathbb{Q}}[\mathcal{K}_i],\:i=1,\ldots,\aleph\right\},
\end{equation}
denote the space of all possible configurations of an McQN with $\aleph$ stations and $d$ classes, in which classes in $\mathcal{K}_i$  are assigned to queue/station $i$; obviously $\mathbb{X}$ is denumerable. Furthermore, we denote the empty configuration on $\mathbb{X}$ by $\emptyset:=[\emptyset,\ldots,\emptyset]$. Finally, for $\boldsymbol{\xi}=[p_1,\ldots,p_\aleph]\in\mathbb{X}$ we define $\mathcal{V}_k(\boldsymbol{\xi}):=\mathcal{V}_k(p_i)$, provided that $k\in\mathcal{K}_i$.

On the class of real-valued functions $\mathbb{R}^\mathbb{X}$ we define the following linear operators:
\begin{itemize}
\item for $h:\mathbb{X}\longrightarrow\mathbb{R}$ we define the $h$-multiplication operator $\Psi[h]\phi:=h\cdot\phi$;
\item for $f:\mathbb{X}\longrightarrow\mathbb{X}$ we define the $f$-composition operator $\Phi[f]\phi:=\phi\circ f$.
\end{itemize}

We are now in the position to formally define the concept of Q-process.

\begin{defn}\label{Q:def} \emph{Consider the following numerical elements:
\begin{itemize}
\item two vectors $\boldsymbol{\theta}=(\theta_k:k=1,\ldots,d)\geq\mathbf{0}$, $\boldsymbol{\beta}=(\beta_k:k=1,\ldots,d)>\mathbf{0}$.
\item a {sub-stochastic matrix} $\boldsymbol{R}=\{R_{kl}\}_{k,l=1,\ldots,d}$, satisfying (\ref{R-series:eq}).
\end{itemize}
A Q-process defined by parameter  $(\boldsymbol{\theta},\boldsymbol{\beta},\boldsymbol{R})$ is the continuous-time Markov chain $\mathcal{X}:=\{X_t\}_{t\geq 0}$
on the space $\mathbb{X}$, given by (\ref{state-sp:eq}), having generator
\begin{equation}\label{generator:eq}
{\boldsymbol A}:=\sum_{(k,l)\in\mathcal{T}}\Psi\left[h_{(k,l)}\right]\left(\Phi[f_{(k,l)}]-{\boldsymbol I}\right),
\end{equation}
where $\mathcal{T}:=\{0,1,\ldots,d\}^2\setminus\{(0,0)\}$ and for $\boldsymbol{\xi}=[p_1,\ldots,p_\aleph]\in\mathbb{X}$ we define:
\begin{itemize}
\item $f_{(0,k)}(\boldsymbol{\xi})=[p_1,\ldots,\mathcal{I}_k(p_i),\ldots,p_\aleph]$ and $f_{(k,0)}(\boldsymbol{\xi})=[p_1,\ldots,\mathcal{D}_k(p_i),\ldots,p_\aleph]$;
\item $f_{(k,l)}=f_{(0,l)}\circ f_{(k,0)}$, for $l\neq 0$;
\item $h_{(0,k)}(\boldsymbol{\xi})=\theta_k$ and $h_{(k,l)}(\boldsymbol{\xi})=\beta_k \mathcal{V}_k(p_i) R_{kl}$,
\end{itemize}
for $k\in\mathcal{K}_i$, $i=1,\ldots,\aleph$ and $l=0,1,\ldots,d$.}\hfill$\diamond$
\end{defn}

Regarding Definition \ref{Q:def}, a few remarks are in order.
\begin{itemize}
\item In the above definition, $(0,k)$-transitions correspond to external arrivals to class $k$, $(k,0)$-transitions correspond
to external departures from class $k$, while $(k,l)$-transitions correspond to switches from class $k$ to class $l$.
\item For any $(k,l)\in\mathcal{T}$, $h_{(k,l)}(\boldsymbol{\xi})$ gives the transition rate corresponding to a $(k,l)$-transition from state
$\boldsymbol{\xi}$, while $f_{(k,l)}(\boldsymbol{\xi})$ denotes the state-space transform after performing a $(k,l)$-transition from $\boldsymbol{\xi}$. 
\item In some situations, depending on the individual service disciplines, the resulting Q-process is lumpable, i.e., the state-space
$\mathbb{X}$ can be partitioned in equivalence classes and the resulting quotient process is still Markov. In such cases, the reduced models agree with the ones in \cite{Dai:95}; see Appendix \ref{space:sec}.  
\item A Q-process is called \emph{elementary} if $\aleph=d$. Elementary Q-processes correspond to Jackson network models. They are lumpable with quotient space $\mathbb{N}^d$.
\end{itemize}

In accordance with the CTMC formalism, we shall introduce the following notations: for a Q-process $\mathcal{X}=\{X_t:t\geq 0\}$
on $\mathbb{X}$ we denote by $\boldsymbol{P}^t$ the associated transition operator, defined as
(for suitable $\phi:\mathbb{X}\longrightarrow\mathbb{R}$)
$$\forall t\geq 0,\boldsymbol{\xi}\in\mathbb{X}:\:\boldsymbol{P}^t(\boldsymbol{\xi},\phi):=\mathbb{E}[\phi(X_t)|X_0=\boldsymbol{\xi}]=\mathbb{E}^{\boldsymbol{\xi}}[\phi(X_t)].$$
We further set $\boldsymbol{P}^t(\boldsymbol{\xi},\Omega):=\boldsymbol{P}^t(\boldsymbol{\xi},\mathbf{1}_\Omega)$, for $\Omega\subseteq\mathbb{X}$. The reader should bear in mind that
a Q-process depends (numerically) on the parameters $(\boldsymbol{\theta},\boldsymbol{\beta},\boldsymbol{R})$; {to emphasize the dependence on $\boldsymbol{\theta}$, we use the notation $\boldsymbol{P}_{\boldsymbol{\theta}}^t$ and $\mathbb{E}_{\boldsymbol{\theta}}$}.

Note that the generator in (\ref{generator:eq}) defines a bounded linear operator on $\mathcal{C}_0[\mathbb{X}]$, i.e.,  the space
of functions $\phi\in\mathbb{R}^\mathbb{X}$ vanishing at infinity\footnote{$\phi:\mathbb{X}\longrightarrow\mathbb{R}$ is vanishing
at infinity if there exists an increasing sequence of exhausting compacts $\{\Omega_n\}_{n\in\mathbb{N}}$ satisfying $\sup\{|\phi(\boldsymbol{\xi})|:\boldsymbol{\xi}\notin\Omega_n\}\longrightarrow 0$, for $n\rightarrow\infty$.}, having norm 
$$\|{\boldsymbol A}\|=\sum_{k=1}^d\theta_k+\sum_{i=1}^\aleph\max_{k\in\mathcal{K}_i}\beta_k<\infty.$$
In particular, for $\phi\in\mathcal{C}_0[\mathbb{X}]$ it holds that 
\begin{equation}\label{semigroup:eq}
\forall t\geq 0,\boldsymbol{\xi}\in\mathbb{X}:\:\boldsymbol{P}^t(\boldsymbol{\xi},\phi)=[\exp(t{\boldsymbol A})\phi](\boldsymbol{\xi}),
\end{equation}  

\begin{remk}\emph{Since $\|{\boldsymbol A}\|<\infty$ the (bounded) linear operator ${\boldsymbol A}$ extends in a natural way
to the class of bounded functions in $\mathbb{R}^\mathbb{X}$ (having the same norm) and so does $\exp(t{\boldsymbol A})$, hence
(\ref{semigroup:eq}) extends to all bounded $\phi$'s.}\hfill$\diamond$
\end{remk}

Furthermore, ${\boldsymbol Q}:={\boldsymbol I}+(1/a){\boldsymbol A}$ defines a Markov (transition) operator on the class
of bounded functions on $\mathbb{X}$, for any $a\geq\|{\boldsymbol A}\|$ and
\begin{equation}\label{embed:prob}
\boldsymbol{P}^t=\exp(t{\boldsymbol A})=\exp[a t({\boldsymbol Q}-{\boldsymbol I})]=
\exp(-a t)\sum_{n\geq 0}\frac{(a t)^n}{n!}{\boldsymbol Q}^n.
\end{equation}
The above property is called \emph{uniformization} and allows one to sample the process $\mathcal{X}$ via the so-called {uniformized} (Markov) $a$-\emph{chain}, as follows: letting $\Xi=\{\Xi_n:n\geq 0\}$ denote the Markov chain
with transition operator $\boldsymbol{Q}$, a random sample from $X_t$ can be obtained as $\Xi_{N_t}$, where $\{N_t:t\geq 0\}$ denotes a Poisson process with rate $a$; that is, $X_t$ and $\Xi_{N_t}$ coincide in distribution, as readily follows from (\ref{embed:prob}). 


\section{Stochastic Monotonicity of Q-processes}\label{mon:sec}
In this section we first introduce a (stochastic) monotonicity concept tailored to Q-processes and then deduce further properties
of Q-processes satisfying such monotonicity assumptions. In addition, we show that this type of monotonicity is quite common for
Q-processes, by identifying a pair of regularity conditions on the underlying queue policies and service allocations  which guarantee monotonicity of a Q-process. Importantly, these conditions are met by virtually all usual queue policies and service allocations used
in applications; this is shown in Appendix \ref{space:sec}.\vspace{2mm}

If $\mathcal{K}$ is an arbitrary set of classes, we define the canonical partial ordering $\subseteq$ on $\overline{\mathbb{Q}}[\mathcal{K}]$ as follows: $p\subseteq q$ if the digits of $p$ can be identified among the digits of $q$, \emph{in the same order}. Formally, $\emptyset\subseteq p$, for any $p$ and if $p=(k_1,\ldots,k_m)$ and $q=(l_1,\ldots,l_n)$, with $1\leq m\leq n$, then $p\subseteq q$ iff there exists some increasing sequence $\nu_1<\ldots<\nu_m$ satisfying $k_\imath=l_{\nu_\imath}$, for $\imath=1,\ldots,m$. In addition, $p$ and $q$ are called \emph{consecutive} if $n=m+1$; in this case, $q$ can be obtained by inserting an extra digit in the sequence $p$, in some arbitrary position. Finally, we note that $\mathcal{D}_k(p)\subseteq p\varsubsetneq\mathcal{I}_l(p)$, for $k,l=1,\ldots,d$, and that $p$ and $\mathcal{I}_l(p)$ are consecutive sequences for any queue
policy. 

Furthermore, we extend $\subseteq$ to $\mathbb{X}$, defined in (\ref{state-sp:eq}), as follows: if $\boldsymbol{\xi}=[p_1,\ldots,p_\aleph]\in\mathbb{X}$
and $\boldsymbol{\zeta}=[q_1,\ldots,q_\aleph]\in\mathbb{X}$ then $\boldsymbol{\xi}\subseteq\boldsymbol{\zeta}$ iff $p_i\subseteq q_i$, for $i=1,\ldots,\aleph$. Also, we say that $(\boldsymbol{\xi},\boldsymbol{\zeta})\in\mathbb{X}\times\mathbb{X}$ is a pair of \emph{consecutive configurations} if there exists some $i_0$ such that $p_i=q_i$ for $i\neq i_0$, whereas for $i=i_0$ the sequences $p_i$ and $q_i$ are consecutive. Note that for any such pair of consecutive
configurations there exists exactly one $b\in\mathcal{K}_{i_0}$ such that $p_{i_0}$ and $q_{i_0}$ differ by exactly one $b$-digit; we denote by $\Delta_b\subseteq\mathbb{X}\times\mathbb{X}$ the set of all (pairs of) consecutive configurations on $\mathbb{X}$ differing
by a $b$-digit, so that the family $\{\Delta_b:b=1,\ldots,d\}$ forms a partition of the set of all (pairs of) consecutive configurations.
Note that $f_{(k,0)}(\boldsymbol{\xi})\subseteq\boldsymbol{\xi}\varsubsetneq f_{(0,l)}(\boldsymbol{\xi})$ and $(\boldsymbol{\xi},f_{(0,l)}(\boldsymbol{\xi}))\in\Delta_l$ for $\boldsymbol{\xi}\in\mathbb{X}$ and $k,l=1,\ldots,d$. Note that, for any $\boldsymbol{\xi}\varsubsetneq\boldsymbol{\zeta}$, there exists a sequence of consecutive configurations $\boldsymbol{\xi}_0,\boldsymbol{\xi}_1,\ldots,\boldsymbol{\xi}_n$, i.e.\ $(\boldsymbol{\xi}_m,\boldsymbol{\xi}_{m+1})$ are consecutive for any $m=0,1,\ldots,n-1$, such that $\boldsymbol{\xi}_0=\boldsymbol{\xi}$ and $\boldsymbol{\xi}_n=\boldsymbol{\zeta}$.

Finally, the mapping $\phi:\mathbb{X}\longrightarrow\mathbb{R}$ is called \emph{increasing} if $\boldsymbol{\xi}\subseteq\boldsymbol{\zeta}$ entails
$\phi(\boldsymbol{\xi})\leq\phi(\boldsymbol{\zeta})$. Note that it suffices to verify that the latter property only holds for consecutive configurations $(\boldsymbol{\xi},\boldsymbol{\zeta})$. In particular, if $\phi$ is increasing then $(\Phi[f_{(0,k)}]-{\boldsymbol I})\phi\geq\mathbf{0}$ (pointwise), for any $k=1,\ldots,d$. 

In what follows, we denote by $\mathcal{J}[\mathbb{X}]$ the class of increasing functions on $\mathbb{X}$ and let $\mathcal{F}\subseteq\mathcal{J}[\mathbb{X}]$. Our next definition introduces the concept of $\mathcal{F}$-monotonicity.

\begin{defn}\label{mon:def}
\emph{The Q-process $\mathcal{X}$ having generator ${\boldsymbol A}$ is called $\mathcal{F}$-monotone if there exists some $a\geq\|\boldsymbol{A}\|$ such that the transition operator
$${\boldsymbol Q}^n=[{\boldsymbol I}+(1/a){\boldsymbol A}]^n,$$ maps $\mathcal{F}$ onto $\mathcal{J}[\mathbb{X}]$,
for any $n\geq 0$; more specifically, ${\boldsymbol A}$ is called $\mathcal{F}$-monotone if $\boldsymbol{\xi}\subseteq\boldsymbol{\zeta}$ and
$\phi\in\mathcal{F}$ entails $[{\boldsymbol Q}^n\phi](\boldsymbol{\xi})\leq[{\boldsymbol Q}^n\phi](\boldsymbol{\zeta})$, for any $n\geq 0$}.\hfill$\diamond$
\end{defn}

If $\mathcal{F}=\mathcal{J}[\mathbb{X}]$ in Definition \ref{mon:def} then we call $\mathcal{X}$ \emph{strongly} monotone.

\begin{remk} 
\emph{$\mathcal{F}$-monotonicity of a Q-process requires that the transition operator ${\boldsymbol Q}^n$ maps $\mathcal{F}$ onto $\mathcal{J}[\mathbb{X}]$, for $n\geq 0$. In particular, if the statement in Definition \ref{mon:def} holds true for some $a \geq\|\boldsymbol{A}\|$ then it holds for any $a'\geq a$, since
$$\boldsymbol{I}+\frac{1}{a'}\boldsymbol{A}=\frac{(a'-a)}{a'}\boldsymbol{I}+\frac{a}{a'}\boldsymbol{Q}.$$
Furthermore, we note that $\mathcal{F}$-monotonicity of a Q-process $\mathcal{X}=\{X_t:t\geq 0\}$ entails 
$$\forall t\geq 0:\:\boldsymbol{P}^t(\boldsymbol{\xi},\phi)\leq \boldsymbol{P}^t(\boldsymbol{\zeta},\phi),$$ for any $\boldsymbol{\xi}\subseteq\boldsymbol{\zeta}$ and $\phi\in\mathcal{F}$; this readily follows
from (\ref{semigroup:eq}) and (\ref{embed:prob})}.\hfill$\diamond$
\end{remk}

$\mathcal{F}$-monotonicity is essentially different (in fact, weaker) from similar concepts introduced in \cite{Massey:87}, where stochastic monotonicity amounts to $\boldsymbol{Q}$ leaving invariant some $\mathcal{F}\subseteq\mathcal{J}[\mathbb{X}]$. In fact, strong monotonicity (in the sense of Definition \ref{mon:def}) is equivalent to the strong monotonicity concept introduced in \cite{Massey:87}. However, we stress that, in general, for $\mathcal{F}\neq\mathcal{J}{\mathbb{X}}$, the monotonicity concept in Definition \ref{mon:def} is weaker than that introduced in \cite{Massey:87}. 

Although weaker than the standard stochastic monotonicity, $\mathcal{F}$-monotonicity still has rather powerful implications, as shown by our next result.

\begin{prop}\label{mon:prop}
\emph{Let $\mathcal{X}$ denote a Q-process which is $\mathcal{F}$-monotone, for some given parameter $(\boldsymbol{\theta},\boldsymbol{\beta},\boldsymbol{R})$, for some $\mathcal{F}\subseteq\mathcal{J}[\mathbb{X}]$. Then:}
\begin{enumerate}
\item [{\rm I}.] \emph{For every $\phi\in\mathcal{F}$, the mapping $t\longmapsto \boldsymbol{P}_{\boldsymbol{\theta}}^t(\emptyset,\phi)$ is non-decreasing.}
\item [{\rm II}.] \emph{If $\boldsymbol{\theta}\leq\boldsymbol{\vartheta}$ (resp.\ $\geq$) then $\boldsymbol{P}_{\boldsymbol{\theta}}^t(\boldsymbol{\xi},\phi)\leq\boldsymbol{P}_{\boldsymbol{\vartheta}}^t(\boldsymbol{\xi},\phi)$ (resp.\ $\geq$), for every $\boldsymbol{\xi}\in\mathbb{X},\phi\in\mathcal{F}$.}

\noindent \emph{In particular, if $\mathcal{X}$ is $\mathcal{F}$-monotone for any $\boldsymbol{\theta}$, the mapping $\boldsymbol{\theta}\longmapsto\boldsymbol{P}_{\boldsymbol{\theta}}^t(\boldsymbol{\xi},\phi)$ is non-decreasing, for every $\boldsymbol{\xi}\in\mathbb{X},\phi\in\mathcal{F}$.\hfill$\diamond$}
\end{enumerate}
\end{prop}

\subsection{\texorpdfstring{$\mathcal{F}$}{F}-\textsf{\emph{Monotonicity for Elementary Q-processes (Jackson Networks)}}} Note that, for elementary Q-processes all service disciplines are equivalent (resulting in the same process -- see Appendix \ref{space:sec}), so that the state-space 
$\mathbb{X}$ can be reduced to $\mathbb{N}^d$ and $\subseteq$ corresponds to the usual component-wise ordering $\leq$ on $\mathbb{N}^d$.

We claim that $\boldsymbol{Q}=\boldsymbol{I}+(1/a)\boldsymbol{A}$ maps $\mathcal{J}[\mathbb{N}^d]$ onto itself, for any $a\geq\|\boldsymbol{A}\|$, for any parameter $(\boldsymbol{\theta},\boldsymbol{\beta},\boldsymbol{R})$. In fact, it suffices to prove that
$$\Psi[h_{(k,l)}]\left(\Phi[f_{(k,l)}]-{\boldsymbol I}\right),$$ does that, for each $(k,l)\in\mathcal{T}$.
The last claim follows by  \cite[Theorem 5.4]{Massey:87} by noting that ${\boldsymbol x}\leq{\boldsymbol z}$ entails either
\begin{itemize}
\item $h_{(k,l)}({\boldsymbol x})=h_{(k,l)}({\boldsymbol z})$, in which case it holds that $f_{(k,l)}({\boldsymbol x})\leq f_{(k,l)}({\boldsymbol z})$, or
\item $0=h_{(k,l)}({\boldsymbol x})<h_{(k,l)}({\boldsymbol z})$, in which case $k\neq 0$ and $x_k=0$, hence 
$$\forall l=0,1,\ldots,d:\:{\boldsymbol x}\leq f_{(k,0)}({\boldsymbol z})\leq f_{(k,l)}({\boldsymbol z}).$$
\end{itemize}
One concludes that elementary Q-processes are strongly monotone, implying that ${\boldsymbol x}\leq{\boldsymbol z}$ entails 
$\boldsymbol{P}^t({\boldsymbol x},\phi)\leq \boldsymbol{P}^t({\boldsymbol z},\phi)$, for any $t\geq 0$ and $\phi\in\mathcal{J}[\mathbb{N}^d]$.

\newpage

\subsection{\texorpdfstring{$\mathcal{F}$}{F}-\textsf{\emph{Monotonicity for Non-Elementary Q-processes}}} \label{mc-mon:sec}

In this subsection we discuss $\mathcal{F}$-monotonicity in the non-elementary framework, i.e.,  for Q-processes corresponding
to McQNs where at least one server is multi-class. It turns out that strong monotonicity does not carry over beyond the
elementary framework, the main reason being that the same transition probabilities from a given pair of ordered states are
\emph{not} stochastically ordered (with respect to any sensible stochastic ordering which is compatible with $\subseteq$
on $\mathbb{X}$), hence the results in \cite{Massey:87} do not apply in this framework.

We shall establish, instead, a weaker form of monotonicity. More specifically, for some arbitrary set of classes $\mathcal{K}$ let
us denote by $\ell:\mathbb{Q}[\mathcal{K}]\longrightarrow\mathbb{N}$ the mapping assigning to the sequence $p:=(k_1,\ldots,k_n)$
its length $\ell(p):=n$. We extend the length mapping, as follows: first we extend it to $\overline{\mathbb{Q}}[\mathcal{K}]$, by 
setting $\ell(\emptyset):=0$ and finally on the state-space $\mathbb{X}$, in (\ref{state-sp:eq}), we define $\ell:\mathbb{X}\longrightarrow\mathbb{N}$ via $$\ell(\boldsymbol{\xi}):=\sum_{i=1}^\aleph\ell(p_i).$$ 
Furthermore, let us consider the space
\begin{equation}\label{G-space:eq}
\mathcal{G}:=\{\phi=G\circ\ell:\:G\in\mathbb{R}^\mathbb{N}\:\text{is non-decreasing}\}\subseteq\mathcal{J}[\mathbb{X}],
\end{equation}
i.e., the class of functions that are non-decreasing in the cumulated queue length.

Our next result establishes sufficient conditions for $\mathcal{G}$-monotonicity.

\begin{them}\label{HQ:them}
\emph{Let $\mathcal{X}$ denote a Q-process (cf. Definition \ref{Q:def}) satisfying:
\begin{enumerate}
\item [(i)] for any $\boldsymbol{\xi}\subseteq\boldsymbol{\zeta}$ and $k=1,\ldots,d$ it holds $f_{(0,k)}(\boldsymbol{\xi})\subseteq f_{(0,k)}(\boldsymbol{\zeta})$.
\item [(ii)] for any pair of consecutive configurations $(\boldsymbol{\xi},\boldsymbol{\zeta})\in\Delta_b$, with $b=1,\ldots,d$, it holds that
$\mathcal{V}_k(\boldsymbol{\zeta})\leq \mathcal{V}_k(\boldsymbol{\xi})$, for $k\neq b$; in particular, $\mathcal{V}_b(\boldsymbol{\xi})\leq \mathcal{V}_b(\boldsymbol{\zeta})$.
\end{enumerate} 
Then $\mathcal{X}$ is $\mathcal{G}$-monotone, with $\mathcal{G}$ defined in (\ref{G-space:eq}), for any parameter $(\boldsymbol{\theta},\boldsymbol{\beta},\boldsymbol{R})$}.\hfill$\diamond$
\end{them}

Conditions (i) and (ii) in Theorem \ref{HQ:them} impose some monotonicity assumptions on the insertion operators, resp.\ on
the service allocations, associated with the Q-process and can be verified for each server individually. Condition~(i) requires
that new arrivals do not influence the order of the jobs (already) in the queue (i.e., the insertion operators are order-preserving), whereas Condition (ii) requires that the service capacity allocated to any class may not decrease by increasing the number of its representatives. Such conditions are fulfilled by a wide range of service disciplines used in queueing applications, e.g. first-came-first-served, static buffer priority or standard processor sharing disciplines; see Appendix~\ref{space:sec} for details.
As such, Theorem \ref{HQ:them} covers virtually all McQNs of practical interest.

Note that, if $b\in\mathcal{K}_{i_0}$ then $\mathcal{V}_k(\boldsymbol{\xi})=\mathcal{V}_k(\boldsymbol{\zeta})$ for any $k\notin\mathcal{K}_{i_0}$, since the two configurations are
identical at all servers $i\neq i_0$; therefore, the inequalities in Condition (ii) are only relevant for $k\in\mathcal{K}_{i_0}$. 

\section{Stability of \texorpdfstring{$\mathcal{F}$}{F}-Monotone Q-processes}\label{stab:sec}

In this section we discuss stability of Q-processes and its relation to $\mathcal{F}$-monotonicity; here stability is understood in a Markovian sense. In particular, we consider a Q-process $\mathcal{X}$ with parameter  $(\boldsymbol{\theta},\boldsymbol{\beta},\boldsymbol{R})$ and regard it as a parametric Markov process, controlled by the arrival rate vector $\boldsymbol{\theta}\in\Theta$; hence $\boldsymbol{\beta}$ and $\boldsymbol{R}$ are kept fixed. For such a process we show that stability is a monotone property with respect to $\boldsymbol{\theta}$.

To start with, we note that a Q-process is irreducible iff $(\boldsymbol{I}-\boldsymbol{R}')^{-1}\boldsymbol{\theta}>\mathbf{0}$ and the irreducibility support, denoted by $\mathbb{X}_0$, depends on the underlying queue policies, but not on $\boldsymbol{\theta}$. In what follows, we let $\Theta\subseteq\{\boldsymbol{\theta}\geq\mathbf{0}:(\boldsymbol{I}-\boldsymbol{R}')^{-1}\boldsymbol{\theta}>\mathbf{0}\}$, hence assuming that $\mathcal{X}$ is irreducible (on some $\mathbb{X}_0\subseteq\mathbb{X}$), for any $\boldsymbol{\theta}\in\Theta$.

\begin{remk}
\emph{The empty configuration $\emptyset$ always belongs to the irreducibility support $\mathbb{X}_0$, since $\emptyset$ is attainable regardless of the queue policies}.\hfill$\diamond$
\end{remk}

We know from standard theory that any irreducible continuous-time Markov chain is either transient or (null or positive) 
recurrent  \cite[{Thm.\ 8.2.5}]{MT:93}; in addition, positive recurrence is equivalent to stability/ergodicity
 \cite[{Thm.\ 13.3.3}]{MT:93} and guarantees the existence of an (essentially) unique equilibrium (limiting) distribution supported
on $\mathbb{X}_0$. We define the \emph{stability region} associated with a Q-process $\mathcal{X}$ by
$$\Theta_{\rm s}:=\{\boldsymbol{\theta}\in\Theta:\:\mathcal{X}\:
\text{is positive recurrent under}\:\mathbb{P}_{\boldsymbol{\theta}}\};$$
here $\mathbb{P}_{\boldsymbol{\theta}}$ denotes the probability law of the Q-process. An alternative way to characterize stability is as
follows: let $T_{\boldsymbol{\xi}}:=\inf\{t>0:\:X_t=\boldsymbol{\xi},\:X_{t-}\neq\boldsymbol{\xi}\}$, for $\boldsymbol{\xi}\in\mathbb{X}$ denote the (first) hitting time of the state $\boldsymbol{\xi}$
({after the process visited at least one different state}). Then the process $\mathcal{X}$ is stable if and only if $\mathbb{E}_{\boldsymbol{\theta}}^{\boldsymbol{\xi}}[T_{\boldsymbol{\xi}}]<\infty$, for any $\boldsymbol{\xi}\in\mathbb{X}_0$. Moreover, it suffices that the expected return time is finite (only) for some particular $\boldsymbol{\xi}$; see e.g. \cite{Meyn-Tweedie:II}.


For any bounded function $\phi:\mathbb{X}\longrightarrow\mathbb{R}$ and measure $\mu$ on $\mathbb{X}$ we set
\begin{equation}\label{cycle:eq}
\langle\mu,\phi\rangle:=\int\phi(\boldsymbol{\xi})\mu(d\boldsymbol{\xi}),\quad\mathcal{R}_{\boldsymbol{\theta}}[\phi]:=
\mathbb{E}_{\boldsymbol{\theta}}^\emptyset\left[\int_0^{T_\emptyset}\phi(X_t)\:dt\right].
\end{equation}
Since $T_{\boldsymbol{\xi}}$ is always larger than the holding time of $\mathcal{X}$ in the initial state $X_0$, we have
{$\mathcal{R}_{\boldsymbol{\theta}}[\mathbf{1}]=\mathbb{E}_{\boldsymbol{\theta}}^\emptyset[T_\emptyset]\geq 1/\|\boldsymbol{\theta}\|$}, hence stability amounts to
$\mathcal{R}_{\boldsymbol{\theta}}[\mathbf{1}]<\infty$. Furthermore, for $\Theta\in\boldsymbol{\theta}_{\rm s}$ and any bounded $\phi$ it holds that $\mathcal{R}_{\boldsymbol{\theta}}[\phi]<\infty$
and $\mathcal{R}_{\boldsymbol{\theta}}[\phi]$ is continuous (even differentiable) in $\boldsymbol{\theta}$.
In addition, denoting by $\pi_{\boldsymbol{\theta}}$ the limiting (equilibrium) distribution of $\mathcal{X}$ under $\mathbb{P}_{\boldsymbol{\theta}}$, the
regenerative ratio formula $\langle\pi_{\boldsymbol{\theta}},\phi\rangle=\mathcal{R}_{\boldsymbol{\theta}}[\phi]/\mathcal{R}_{\boldsymbol{\theta}}[\mathbf{1}]$ holds; see \cite[p.\ 106]{AG:2007}. 
Hence $\boldsymbol{\theta}\longmapsto\langle\pi_{\boldsymbol{\theta}},\phi\rangle$ is continuous on $\Theta_{\rm s}$, for any bounded $\phi\in\mathbb{R}^\mathbb{X}$.
Finally, we note that the stability region $\Theta_{\rm s}$ is an open subset of $\Theta$, since 
$$\Theta_{\rm s}=\bigcup_{r\geq 0}
\{\boldsymbol{\theta}\in\Theta:\mathcal{R}_{\boldsymbol{\theta}}[\:\mathbf{1}\:]<r\}.$$

The main result of this section provides a characterization of the stability region $\Theta_{\rm s}$ of a  $\mathcal{F}$-monotone Q-process.  

\begin{them}\label{regio:them}
\emph{Let $\mathcal{X}=\{X_t:t\geq 0\}$ be a  $\mathcal{F}$-monotone Q-process for all $\boldsymbol{\theta}\in\Theta$ and assume that there exists $\phi:\mathbb{X}\longrightarrow[0,1]$, vanishing at infinity, such that $\phi(\emptyset)\neq 0$ and $-\phi\in\mathcal{F}$. Define $\varphi_t(\boldsymbol{\theta}):=\boldsymbol{P}_{\boldsymbol{\theta}}^t(\emptyset,\phi)=
\mathbb{E}_{\boldsymbol{\theta}}^\emptyset[\phi(X_t)]$, for $t\geq 0$ and $\boldsymbol{\theta}\in\Theta$.}

\emph{Then the family of functions $\varphi_t:\boldsymbol{\theta}\longrightarrow[0,1]$ satisfies:
\begin{enumerate}
\item [{\rm I}.] $(t,\boldsymbol{\theta})\longmapsto\varphi_t(\boldsymbol{\theta})$ is non-increasing in both $t$ and $\boldsymbol{\theta}$ (componentwise);
\item [{\rm II}.] the limit $\varphi:=\lim_{t\rightarrow\infty}\varphi_t$ is continuous and non-increasing, satisfying
\begin{equation}\label{regio:eq}
\Theta_{\rm s}=\{\boldsymbol{\theta}\in\Theta:\varphi(\boldsymbol{\theta})>0\}.
\end{equation}
In particular, $\Theta_{\rm s}\subseteq\Theta$ is open, star-shaped and $\varphi(\boldsymbol{\theta})=\langle\pi_{\boldsymbol{\theta}},\phi\rangle$,
for $\boldsymbol{\theta}\in\Theta_{\rm s}$.
\item [{\rm III}.] The mapping $\varphi=\langle\pi_\cdot,\phi\rangle:\Theta_{\rm s}\longrightarrow(0,1]$ is strictly decreasing.
\hfill$\diamond$
\end{enumerate}}
\end{them}

Theorem \ref{regio:them} provides a characterization of the stability region of a $\mathcal{F}$-monotone Q-process as the support of some continuous, non-increasing functional $\varphi$. In particular, this shows that $\Theta_{\rm s}$ is an open, star-shaped domain having the origin as a vantage point,
i.e., $\boldsymbol{\theta}\in\Theta_{\rm s}$ entails $c \boldsymbol{\theta}\in\Theta_{\rm s}$, for any $c \in(0,1]$. The property is known in the literature as \emph{star-convexity} (a weaker form of convexity), or \emph{monotonicity} \cite{Dai:99}.

\begin{remk} 
\emph{Note that, in Theorem \ref{regio:them}, $\varphi$ is defined by means of some function $\phi$, whereas the stability region is
a fixed set (depending on the process itself, but not on $\phi$). The representation in (\ref{regio:eq}) is valid for any functional $\varphi$ (hence function $\phi$) satisfying the conditions of the theorem. In particular, $\mathcal{G}$-monotone Q-processes (e.g.\ satisfying the assumptions of Theorem \ref{HQ:them}) satisfy the conditions of Theorem \ref{regio:them}; any function $\phi=G\circ\ell$, with $G$ a non-increasing function, vanishing at infinity, e.g.\ $G(x)=\exp(-\alpha x)$, $G(x)=(1+x)^{-\alpha}$ or
$G(x)=\mathbf{1}\{x\leq\alpha\}$ ($\alpha>0$) can be used}.\hfill$\diamond$
\end{remk}

\section{Concluding Remarks and Discussion}\label{sec:concl}

In this paper, we introduced a general Markovian model (a Q-process) for modeling the dynamics of a Markovian McQN over time. In addition, we introduced a new concept of stochastic ($\mathcal{F}$-) monotonicity. Then we proved that this property holds for a wide class of Q-processes, corresponding to virtually all McQN models which are relevant in applications, e.g.\ McQNs in which any server employs either a first-come-first-served, static buffer priority or processor sharing service discipline. Furthermore, we proved that $\mathcal{F}$-monotonicity is strong enough to ensure that stability is a monotonic property with respect to  the external arrival rates.

The key result of this paper is Theorem \ref{HQ:them} which shows that, under some mild monotonicity assumptions over the insertion operators and service allocations, the configuration-length (total job population in the network) is stochastically monotone with respect to arrival rates and with time (the latter holds provided that the process starts in the empty configuration). However, analyzing the proof of Theorem \ref{HQ:them}, we see that stochastic monotonicity holds, in principle, for any increasing functional of the process which \emph{can only} increase through arrival events. The configuration-length function obviously satisfies this condition for any Q-process, but in particular cases it could also be the case for other functionals, e.g.\ the conditional expected workload in the network, given the network configuration. 
 
The results of this paper facilitate the numerical methods developed in \cite{LM:2016} for determining the stability region of an McQN. To be more precise, let $$\boldsymbol{\theta}=\{\boldsymbol{\theta}=r\cdot\vec{v}:r>0\},$$ where $\vec{v}\geq\mathbf{0}$ denotes a $d$-dimensional vector, satisfying $\|\vec{v}\|=1$; that is, $\Theta$ is a one-dimensional manifold (positive direction) in $\mathbb{R}^d$, endowed with the natural ordering. By Theorem \ref{regio:them}, $\Theta_{\rm s}=(\mathbf{0},\boldsymbol{\theta}_*)$, where 
$$\boldsymbol{\theta}_*=\sup\{\boldsymbol{\theta}\in\Theta:\varphi(\boldsymbol{\theta})>0\}=\min\{\boldsymbol{\theta}\in\Theta:\varphi(\boldsymbol{\theta})=0\};$$ the value $r_*:=\|\boldsymbol{\theta}_*\|$ is called the \emph{stability threshold along direction} $\vec{v}$. Furthermore, the mapping $\varphi:(\mathbf{0},\boldsymbol{\theta}_*)\longrightarrow[0,1]$, $\varphi(\boldsymbol{\theta})=\langle\pi_{\boldsymbol{\theta}},\phi\rangle$, is 
strictly decreasing, hence for any $\varepsilon\in(0,1)$ there exists some unique $\boldsymbol{\theta}_\varepsilon\in(\mathbf{0},\boldsymbol{\theta}_*)$
satisfying $\varphi(\boldsymbol{\theta}_\varepsilon)=\varepsilon$; the value $r_\varepsilon:=\|\boldsymbol{\theta}_\varepsilon\|$ is called the $\varepsilon$-\emph{congestion
threshold along direction} $\vec{v}$. Finally, it is immediate that $r_\varepsilon\uparrow r_*$, as $\varepsilon\downarrow 0$, hence the $\varepsilon$-congestion threshold approximates the stability threshold for $\varepsilon\downarrow 0$.

Numerical (simulation-based) methods for evaluating congestion (and stability) thresholds were developed in \cite{LM:2016} under
some minimal monotonicity assumptions, where stochastic approximation schemes of Robbins-Monro type were applied to some specific
McQN examples. The results in this paper formally validate the numerical results in \cite{LM:2016}.

On the other hand, our analysis shows that concepts such as `stochastic monotonicity' carry over from Jackson networks
to the more general multi-class framework, albeit in a weaker form. Here it is stressed that this weaker form still retains
the most interesting monotonicity features. From a practical standpoint the results in this paper  cover
by and large all relevant cases, but at the theoretical level there are still some questions left. For instance, an interesting
question would be whether there exist Q-processes (resp.\ McQNs) which are {\it not} $\mathcal{G}$-monotone; one might
expect that instances of non-monotone Q-processes could possibly be obtained by violating the conditions of Theorem \ref{HQ:them}.
Finally, we note that such monotonicity properties do not extend to non-Markovian McQN models (which involve non-exponential distributions), as shown in \cite{Whitt:93}.

\newpage 

\appendix

\section{On Queue Policies and Service Allocations}\label{space:sec} 

In this appendix we illustrate the formalism introduced in Section \ref{q:sec}, by providing explicit details
on how the most common service disciplines used in queueing applications fit into our modeling paradigm. In addition, we identify usual queue policies and service allocations satisfying the assumptions of Theorem \ref{HQ:them} (hence giving rise to monotone Q-processes) and also show how the model complexity can be reduced in some special cases. 

Recall that $\mathbb{Q}[\mathcal{K}]$ denotes the space of finite sequences over the set $\mathcal{K}$, defined by (\ref{smcc:eq}). An insertion operator $\mathcal{I}_k$ on $\mathbb{Q}[\mathcal{K}]$ is a mapping such that for any $p$ there exists a decomposition $p=(p',p'')$ such that $\mathcal{I}_k(p)=(p',k,p'')$; a family $\{\mathcal{I}_k:k\in\mathcal{K}\}$ is called a queue policy.

Perhaps the most natural queue policy is `first-came-first-served' (FCFS). In our modeling paradigm, this is obtained as follows: for any $p\in\mathbb{Q}[\mathcal{K}]$, the insertion operator $\mathcal{I}_k$ inserts a $k$-digit at the end of the sequence, i.e.\ $\mathcal{I}_k(p)=(p,k)$, corresponding to the trivial decomposition $p'=p$ and $p''=\emptyset$.

An important class of queue policies widely used in applications is the class of priority-based policies. To formalize that, we call a \emph{priority ranking} on $\mathcal{K}$ a partition $\{\mathcal{C}_1,\ldots,\mathcal{C}_\nu\}$ of $\mathcal{K}$;
this induces a natural (partial) ordering on $\mathcal{K}$, as follows: $k\prec l$ iff there exist $1\leq\imath<\jmath\leq\nu$
such that $k\in\mathcal{C}_\imath$ and $l\in\mathcal{C}_\jmath$. The interpretation is that each partition block
$\mathcal{C}_\imath$ represents a \emph{caste} (subset of unranked classes), with representatives of higher castes
(corresponding to smaller indexes) being allowed to overtake (in the queue) representatives of lower castes (larger indexes).
When every $\mathcal{C}_\imath$ is a singleton the priority ranking is called \emph{total}, as $\prec$ becomes a total ordering
on $\mathcal{K}$; on the other hand, for $\nu=1$ the ordering is trivial. 

We call a \emph{priority policy} a queue policy $\{\mathcal{I}_k:k\in\mathcal{K}\}$ such that for any $p\in\mathbb{Q}[\mathcal{K}]$ there exists some priority ranking (set of castes) so that $\mathcal{I}_k(p)=(p',k,p'')$, where $p''$ is the maximal tail-sequence consisting of (consecutive) digits belonging to lower castes than $k$; that is, $k$ will overtake all digits with lower priority ranking. Within each caste, a FCFS policy applies. When there is only one caste the queue policy reduces to FCFS. Furthermore, if the priority ranking does not depend on $p$, we call it \emph{static}; otherwise, we call it \emph{dynamic}. Finally, if $p''$ may include  the first digit of the sequence then the policy is called \emph{preemptive}; otherwise, we call it \emph{non-preemptive}.

On the other hand, recall that a service allocation is a state-dependent probability vector on $\mathcal{K}$, specifying what
fraction of service capacity is assigned to any class. An important class of service allocations is the class of \emph{head-of-the-queue} (HQ) allocations; that is, if we define $\kappa:\mathbb{Q}[\mathcal{K}]\longrightarrow\mathcal{K}$ as $\kappa(k_1,\ldots,k_n):=k_1$ then the HQ allocation is defined by $\mathcal{V}_k(p)=\mathbf{1}\{\kappa(p)=k\}$. HQ allocations correspond
to service disciplines in which only one job may receive service at a time. A service allocation which is not of HQ type, is called a \emph{bulk service} (BS) allocation. Among the non-HQ (BS) policies, the most important in applications are the so-called \emph{processor sharing} (PS) service allocations, which only depend on $p$ by means of its \emph{composition vector} $\langle p\rangle$, i.e., the vector whose $k$-th component $\langle p\rangle_k$ denotes the number of $k$-digits in the sequence $p$; in particular, the ordering of the sequence is irrelevant.

By definition, for any PS allocation $W$ there exists some $w:\mathbb{N}^\mathcal{K}\setminus\{\mathbf{0}\}\longrightarrow\mathcal{P}[\mathcal{K}]$,
such that $\mathcal{V}(p)=\mathcal{V}(\boldsymbol{x})$, where $\boldsymbol{x}=\langle p\rangle$. Usual choices are:
\begin{itemize}
\item \emph{Egalitarian} allocation, specified by the mapping
$$\mathcal{V}_k({\boldsymbol x})=\frac{\mathbf{1}\{k\in\varsigma[{\boldsymbol x}]\}}{\#\varsigma[{\boldsymbol x}]},$$
i.e., the server capacity is uniformly distributed among classes.
\item \emph{Proportional} allocation, specified by the mapping
$$\mathcal{V}_k({\boldsymbol x})=\frac{x_k}{\|{\boldsymbol x}\|},$$
i.e., the server capacity is distributed proportionally with the number of representatives in each class.
\item \emph{Preferential} allocation (assumes a total priority order on $\mathcal{K}$), given by
$$\mathcal{V}_k({\boldsymbol x})=\mathbf{1}\{\kappa({\boldsymbol x})=k\},$$
where $\kappa({\boldsymbol x})$ denotes the highest-ranked class in $\varsigma({\boldsymbol x})$.
\end{itemize}

The usual FCFS and SBP (static buffer priority) service disciplines are recovered in our model via the respective queue policies, 
in combination with a HQ service allocation. One can easily verify that FCFS and static priority policies (of both preemptive and non-preemptive type) satisfy Condition (i) in Theorem \ref{HQ:them}, whereas dynamic priority policies do not, in general. 
Furthermore, processor sharing disciplines are obtained via the corresponding PS allocation in combination with any queue policy
(which is irrelevant). The PS allocations listed above satisfy Condition~(ii) in Theorem \ref{HQ:them}. One concludes that, indeed,
Markovian models associated with virtually all McQNs of interest in applications are $\mathcal{G}$-monotone, cf.\ Theorem \ref{HQ:them}.

In some situations of interest, it is possible to reduce the complexity of the full space $\mathbb{Q}[\mathcal{K}]$ (and of the full space $\mathbb{X}$, accordingly), by identifying equivalent configurations; this is formalized as follows: the space $(\mathbb{Q}[\mathcal{K}];\mathcal{I}_k,\mathcal{D}_k,\mathcal{V}_k:k\in\mathcal{K})$ of multi-class configurations over the set $\mathcal{K}$ is \emph{reducible} if there exists an equivalence relation $\sim$ on $\mathbb{Q}[\mathcal{K}]$ such that $p\sim q$ entails $\langle p\rangle=\langle q\rangle$ and $\mathcal{V}_k(p)=\mathcal{V}_k(q)$, $\mathcal{I}_k(p)\sim\mathcal{I}_k(q)$ and $\mathcal{D}_k(p)\sim\mathcal{D}_k(q)$, for all $k\in\mathcal{K}$. One can further extend $\sim$ to the augmented space $\overline{\mathbb{Q}}[\mathcal{K}]$ by identifying the empty configuration with itself. The mappings $\mathcal{I}_k,\mathcal{D}_k,\mathcal{V}_k$ are then well defined on the quotient space $\overline{\mathbb{Q}}[\mathcal{K}]/\sim$, which will be called a \emph{reduced space} of multi-class configurations.

Instances of reducible spaces of multi-class configurations are given below:
\begin{enumerate}
\item If $\#\mathcal{K}=1$ (single-class) then $p\sim q$ iff $\langle p\rangle=\langle q\rangle$ gives $\overline{\mathbb{Q}}[\mathcal{K}]/\sim\:=\mathbb{N}$.
In particular, queue policies and service allocations are irrelevant.
\item For a static (total, non-preemptive) priority policy with HQ allocation, let $p\sim q$ iff $\kappa(p)=\kappa(q)$ and $\langle p\rangle = \langle q\rangle$, which gives $\overline{\mathbb{Q}}[\mathcal{K}]/\sim\:=\mathcal{K}\times\mathbb{N}^\mathcal{K}$.
\item For a PS allocation, $p\sim q$ iff $\langle p\rangle=\langle q\rangle$; in this case, $\overline{\mathbb{Q}}[\mathcal{K}]/\sim\:=\mathbb{N}^\mathcal{K}$.
The same factorization holds for static (total, preemptive) priority policies with HQ allocations, which can be recovered via a PS preferential allocation.
\end{enumerate}

\section{Proofs of the Results}

In this appendix we provide the proofs of the results presented in this paper.

\noindent{\bf Proof of Proposition \ref{mon:prop}}: Let ${\boldsymbol A}_{\boldsymbol{\theta}}$, for $\boldsymbol{\theta}\geq\mathbf{0}$, denote the generator  of the Q-process $\mathcal{X}$, defined in (\ref{generator:eq}) with arrival rate vector $\boldsymbol{\theta}$, i.e., with parameter $(\boldsymbol{\theta},\boldsymbol{\beta},\boldsymbol{R})$.

I. By Kolmogorov's Equation, the mapping of interest is differentiable with respect to ~$t$ and, cf.\ (\ref{semigroup:eq}),
it holds that
\begin{equation}\label{Kolm:eq}
\frac{d}{dt}\boldsymbol{P}_{\boldsymbol{\theta}}^t(\emptyset,\phi)=\frac{d}{dt}\:
[\exp(t\boldsymbol{A}_{\boldsymbol{\theta}})\phi](\emptyset)=
[\boldsymbol{A}_{\boldsymbol{\theta}}\exp(t\boldsymbol{A}_{\boldsymbol{\theta}})\phi](\emptyset).
\end{equation}
Since $\mathcal{V}_k(\emptyset)=0$, it follows that 
$$h_{(k,l)}(\emptyset)=\beta_k \mathcal{V}_k(\emptyset)R_{kl}=0,$$ for any $(k,l)$ with $k=1,\ldots, d$, $l=0,1,\ldots,d$, 
hence the r.h.s.\ in (\ref{Kolm:eq}) equals
$$\sum_{k=1}^d
\theta_k\left[(\Phi[f_{(0,k)}]-\boldsymbol{I})\exp(t\boldsymbol{A}_{\boldsymbol{\theta}})\phi\right](\emptyset),$$
which is non-negative, for $t\geq 0$, by assumption; this concludes the first part.

II. Define for $\mathbf{0}\leq\boldsymbol{\theta}\leq\boldsymbol{\vartheta}$ and $0\leq s\leq t$,
$$\boldsymbol{E}(s,t):=\exp[(t-s)\boldsymbol{A}_{\boldsymbol{\vartheta}}]\exp(s\boldsymbol{A}_{\boldsymbol{\theta}}).$$
Since $\boldsymbol{E}(0,t)=\exp(t\boldsymbol{A}_{\boldsymbol{\vartheta}})$ and $\boldsymbol{E}(t,t)=\exp(t\boldsymbol{A}_{\boldsymbol{\theta}})$, one 
obtains
\begin{eqnarray}\nonumber
\exp(t\boldsymbol{A}_{\boldsymbol{\vartheta}})-\exp(t\boldsymbol{A}_{\boldsymbol{\theta}}) & = & \int_0^t-\frac{d}{ds}\:\boldsymbol{E}(s,t)\:ds \\
& = & \label{diff-int:eq}
\int_0^t\exp[(t-s)\boldsymbol{A}_{\boldsymbol{\vartheta}}](\boldsymbol{A}_{\boldsymbol{\vartheta}}-\boldsymbol{A}_{\boldsymbol{\theta}})
\exp(s\boldsymbol{A}_{\boldsymbol{\theta}})\:ds,
\end{eqnarray}
and the claim follows from the fact that (for $s\geq 0$)
\begin{equation}\label{eq:aux}
(\boldsymbol{A}_{\boldsymbol{\vartheta}}-\boldsymbol{A}_{\boldsymbol{\theta}})\exp(s\boldsymbol{A}_{\boldsymbol{\theta}})\phi=
\sum_{k=1}^d(\vartheta_k-\theta_k)(\Phi[f_{(0,k)}]-\boldsymbol{I})\exp(s\boldsymbol{A}_{\boldsymbol{\theta}})\phi,
\end{equation}
with the r.h.s.\ above being nonnegative for any $\phi\in\mathcal{F}$, by assumption. 

The case $\boldsymbol{\theta}\geq\boldsymbol{\vartheta}$ can be treated similarly. This concludes the proof.\hfill$\square$ 

\noindent{\bf Proof of Theorem \ref{HQ:them}}: We prove that $\boldsymbol{Q}=\boldsymbol{I}+(1/a)\boldsymbol{A}$,
with $$a=\|\boldsymbol{\theta}\|+\|\boldsymbol{\beta}\|=\sum_{k=1}^d(\theta_k+\beta_k),$$ satisfies Definition \ref{mon:def}. To this end, let
$\mathcal{E}:=\{A_1,\ldots,A_d,B_1,\ldots,B_d\}$ and consider on $\mathcal{E}$ the probability distribution $\mu$ given by
$$\forall k=1,\ldots,d:\:\mu(A_k):=\frac{\theta_k}{a};\quad\mu(B_k):=\frac{\beta_k}{a}.$$
A random transition of the $a$-chain (as defined in Section \ref{q:sec}) from an arbitrary state $\boldsymbol{\xi}\in\mathbb{X}$ can be constructed as follows: generate a random variable $J$ on $\mathcal{E}$, having distribution $\mu$ and
\begin{enumerate}
\item [(A)] if $J=A_k$ (arrival event) then define $\boldsymbol{\xi}=f_{(0,k)}(\boldsymbol{\xi})$;
\item [(B)] if $J=B_k$ then define $\boldsymbol{\xi}=f_{(k,l)}(\boldsymbol{\xi})$ w.p.\ $\mathcal{V}_k(\boldsymbol{\xi})R_{kl}$ ($l=0,1,\ldots,d$) and let $\boldsymbol{\xi}=\boldsymbol{\xi}$ w.p.\
$1-\mathcal{V}_k(\boldsymbol{\xi})$, independently of $J$.
\end{enumerate}

For an arbitrary sample $\{\Xi_\nu:0\leq\nu\leq n\}$ of the $a$-chain, generated according to (A)--(B) above by means
of a sequence of r.v.'s $J_1,\ldots,J_n$, we consider the sequence of arrival events $\mathcal{A}_\nu[\Xi]=[k_1,\ldots,k_m]$ ($0\leq m\leq\nu\leq n$) observed by the chain in the first $\nu$ steps; that is, one accounts for the arrival events within the sequence of $J$'s and records the underlying classes. Note that the probability of observing a given sequence of arrival events equals $$\mathrm{Pr}\{\mathcal{A}_\nu[\Xi]=[k_1,\ldots,k_m]\}=\binom{\nu}{m}
\frac{\theta_{k_1}\cdot\ldots\cdot\theta_{k_m}\cdot\|\boldsymbol{\beta}\|^{\nu-m}}{a^\nu},$$ 
and does not depend on the initial state of the chain.

Returning to the proof of the theorem, we note that it suffices to construct a pair of $a$-chains
$\{(\Xi_\nu',\Xi_\nu''):0\leq\nu\leq n\}$ satisfying $$(\Xi_0',\Xi_0'')=(\boldsymbol{\xi},\boldsymbol{\zeta}),\:\ell(\Xi_n')\leq\ell(\Xi_n''),\:{\rm a.s.}$$
Since the distribution of the sequence of arrival events does not depend on the initial state (hence it is the same for both chains) it suffices to prove the above statement conditioned on the event that the two chains share the same sequence of arrival events. More specifically, we shall prove that the following statement holds true for any $n\geq 0$: 
\begin{itemize}
\item [$(\mathfrak{H}_n)$] For any pair of initial configurations $\boldsymbol{\xi}\subseteq\boldsymbol{\zeta}$
there exists a pair of $a$-chains $\{(\Xi_\nu',\Xi_\nu''):0\leq\nu\leq n\}$ satisfying 
$$(\Xi_0',\Xi_0'')=(\boldsymbol{\xi},\boldsymbol{\zeta}),\:\mathcal{A}_n[\Xi']=\mathcal{A}_n[\Xi''],\:
\ell(\Xi_n')\leq\ell(\Xi_n''),\:{\rm a.s.}$$ 
\end{itemize}
Note that, in the above statement, one may equivalently assume that either $\boldsymbol{\xi}=\boldsymbol{\zeta}$ or $(\boldsymbol{\xi},\boldsymbol{\zeta})$ is a pair of consecutive 
configurations. More concisely, one may assume that $(\boldsymbol{\xi},\boldsymbol{\zeta})\in\Delta_0\cup\Delta_1\cup\ldots\cup\Delta_d$, where, for convenience
we let $\Delta_0:=\{(\boldsymbol{\xi},\boldsymbol{\xi}):\boldsymbol{\xi}\in\mathbb{X}\}$ denote the diagonal set of $\mathbb{X}$.

We shall prove this claim by induction. The key fact in this proof is that, given $\Xi_0'\subseteq\Xi_0''$, there exists a coupling $(\Xi_1',\Xi_1'')$ satisfying $\Xi_1'\subseteq\Xi_1''$ on the event $\{\Xi_1''\neq\Xi_0''\}$; this is guaranteed by Conditions (i) and (ii). While the technical details of this proof are given below, Figure \ref{proof:fig} displays the reasoning used for proving the induction step.

To begin with, note that the statement $(\mathfrak{H}_0)$ is straightforward. Assume now that $(\mathfrak{H}_n)$ holds true for some
$n\geq 0$. Given $(\boldsymbol{\xi},\boldsymbol{\zeta})\in\Delta_b$, for $b=0,1,\ldots,d$, we shall construct a pair of $a$-chains $\{(\Xi_\nu',\Xi_\nu''):0\leq\nu\leq n+1\}$, satisfying 
\begin{equation}\label{eq:coupling}
(\Xi_0',\Xi_0'')=(\boldsymbol{\xi},\boldsymbol{\zeta}),\:\mathcal{A}_{n+1}[\Xi']=\mathcal{A}_{n+1}[\Xi''],\:
\ell(\Xi_{n+1}')\leq\ell(\Xi_{n+1}''),\:{\rm a.s.}
\end{equation}

For $b=0$ (i.e.\ $\boldsymbol{\xi}=\boldsymbol{\zeta}$) the statement is trivial. Consider now the case $b\neq 0$. The key step in this construction is that, given $\Xi_0'=\boldsymbol{\xi}$ and $\Xi_0''=\boldsymbol{\zeta}$, there exists a coupling $(\Xi_1',\Xi_1'')$ satisfying (recall Condition (ii)):
\begin{itemize}
\item [1.] for $k=1,\ldots,d$ we have $\Xi_1''=f_{(0,k)}(\boldsymbol{\zeta})$ if and only if $\Xi_1'=f_{(0,k)}(\boldsymbol{\xi})$.
\item [2.] for $k\neq b$, $\Xi_1''=f_{(k,l)}(\boldsymbol{\zeta})$ entails $\Xi_1'=f_{(k,l)}(\boldsymbol{\xi})$, for any $l=0,1,\ldots,d$.
\item [3.] $\Xi_1''=f_{(b,l)}(\boldsymbol{\zeta})$ entails either $\Xi_1'=f_{(b,l)}(\boldsymbol{\xi})$ or $\Xi_1'=\boldsymbol{\xi}$, for $l=0,1,\ldots,d$.
\end{itemize}
Therefore, on the event $\{\Xi_1''\neq\Xi_0''\}$, we have $(\Xi_1',\Xi_1'')\in\Delta_b$ (cases 1 and 2 above) and $(\Xi_1',\Xi_1'')\in\Delta_b$ or $(\Xi_1',\Xi_1'')\in\Delta_l$ (in case 3); this follows by Condition (i). On the other hand, given that $\{\Xi_1''=\Xi_0''\}$ (i.e.\ $\Xi''$ does not make an actual jump), it readily follows that $(\Xi_0',\Xi_1'')=(\boldsymbol{\xi},\boldsymbol{\zeta})\in\Delta_b$. As such, we distinguish the following cases:

\begin{enumerate}
\item [I.] If $\Xi_1''\neq\Xi_0''$ then $(\Xi_1',\Xi_1'')\in\Delta_0\cup\Delta_1\cup\ldots\cup\Delta_d$. By $(\mathfrak{H}_n)$ there exists a pair of $a$-chains $\{(\Psi_\nu',\Psi_\nu''):0\leq\nu\leq n\}$ such that
$$(\Psi_0',\Psi_0'')=(\Xi_1',\Xi_1''),\:\mathcal{A}_n[\Psi']=\mathcal{A}_n[\Psi''],\:\ell(\Psi_n')\leq\ell(\Psi_n''),\:{\rm a.s.}$$
Defining further $\Xi_\nu':=\Psi_{\nu-1}'$ and $\Xi_\nu'':=\Psi_{\nu-1}''$, for $\nu=1,\ldots,n+1$, one can easily verify the validity of (\ref{eq:coupling}).
\item [II.] If $\Xi_1''=\Xi_0''$ then $(\Xi_0',\Xi_1'')=(\boldsymbol{\xi},\boldsymbol{\zeta})\in\Delta_b$ and using again $(\mathfrak{H}_n)$ one obtains a pair of $a$-chains $\{(\Psi_\nu',\Psi_\nu''):0\leq\nu\leq n\}$ such that
$$(\Psi_0',\Psi_0'')=(\Xi_0',\Xi_1''),\:\mathcal{A}_n[\Psi']=\mathcal{A}_n[\Psi''],\:\ell(\Psi_n')\leq\ell(\Psi_n''),\:{\rm a.s.}$$
Defining further $\Xi_\nu':=\Psi_\nu'$ and $\Xi_{\nu+1}'':=\Psi_{\nu}''$, for $\nu=1,\ldots,n$, we note that 
$\mathcal{A}_n[\Xi']=\mathcal{A}_{n+1}[\Xi'']$ (since $\mathcal{A}_1[\Xi'']=\emptyset$) and $\ell(\Xi_n')\leq\ell(\Xi_{n+1}'')$.
On the other hand, the constraint $\mathcal{A}_{n+1}[\Xi']=\mathcal{A}_{n+1}[\Xi'']=\mathcal{A}_n[\Xi']$ ensures that the $(n+1)$-st transition of $\Xi'$ does not correspond to an arrival event, whence $$\ell(\Xi_{n+1}')\leq\ell(\Xi_n')\leq\ell(\Xi_{n+1}'');$$ this proves (\ref{eq:coupling}).
\end{enumerate}
Therefore, we proved $(\mathfrak{H}_{n+1})$, which concludes the proof of the theorem.
\hfill$\square$ 

\begin{remk} \emph{Regarding the proof of Theorem \ref{HQ:them}, a few remarks are in order}:
\begin{itemize}
\item \emph{The coupling $\{(\Xi_\nu',\Xi_\nu''):0\leq\nu\leq n\}$ {has the special feature that it} depends on the time-horizon $n$}.
\item \emph{The chain $\Xi''$ is `forced' to perform the same transitions as $\Xi'$, in the same order. Due to the extra initial job, $\Xi''$ may lag behind $\Xi'$, in that the same transitions will be observed later, or is even `dropped' by $\Xi''$; since the two chains are bound to share the same sequence of arrival events, this might result in less (network) departures for $\Xi''$}.
\item \emph{The statement in Theorem \ref{HQ:them} remains valid for any $\phi\in\mathcal{J}[\mathbb{X}]$, satisfying}
$$\forall k\neq 0:\:R_{kl}\cdot(\Phi[f_{(k,l)}]-\boldsymbol{I})\phi\leq 0;$$
\emph{that is, increasing functionals which may only increase through arrival events. Indeed, the key aspect in proving the induction step is that the functional of interest may not increase, given that the chain $\Xi'$ does not see new arrivals}.
\end{itemize}
\end{remk}

\newpage

\noindent{\bf Proof of Theorem \ref{regio:them}}: 
I. It follows directly by Proposition \ref{mon:prop}. 

II. Monotonicity of $\varphi_t(\boldsymbol{\theta})$ with respect to $t$ shows that the limit $\varphi:=\inf_t\varphi_t$ exists and preserves monotonicity with respect to $\boldsymbol{\theta}$. Furthermore, we claim that
\begin{equation}\label{phi:eq}
\varphi(\boldsymbol{\theta})=\left\{
                    \begin{array}{ll}
                      \langle\pi_{\boldsymbol{\theta}},\phi\rangle, & \hbox{$\boldsymbol{\theta}\in\Theta_{\rm s}$;} \\
                      0, & \hbox{$\boldsymbol{\theta}\notin\Theta_{\rm s}$,}
                    \end{array}
                  \right.
\end{equation}
Indeed, since stability entails ergodicity in this context, we have 
$$\forall\boldsymbol{\theta}\in\Theta_{\rm s}:\:\varphi(\boldsymbol{\theta})= \lim_{t\rightarrow\infty}
\mathbb{E}_{\boldsymbol{\theta}}^\emptyset[\phi(X_t)]=\langle\pi_{\boldsymbol{\theta}},\phi\rangle;$$
the r.h.s.\ above is strictly positive since $\phi(\emptyset)>0$ and $\pi_{\boldsymbol{\theta}}$ is supported on $\mathbb{X}_0$, which contains $\emptyset$. On the other hand, let $\boldsymbol{\theta}\notin\Theta_{\rm s}$. Then $P_{\boldsymbol{\theta}}^t(\emptyset,\Omega)\longrightarrow 0$, for every compact (finite) $\Omega\subset\mathbb{X}$; in particular, since $\phi$ is vanishing at infinity, there exists exhausting
compacts $\{\Omega_n\}_{n\in\mathbb{N}}$ such that $\sup\{\phi(\boldsymbol{\xi}):\boldsymbol{\xi}\notin\Omega_n\}\longrightarrow 0$, whence
$$\forall n\in\mathbb{N},t\geq 0:\:
\varphi_t(\boldsymbol{\theta})=\boldsymbol{P}_{\boldsymbol{\theta}}^t(\emptyset,\phi)\leq \boldsymbol{P}_{\boldsymbol{\theta}}^t(\emptyset,\Omega_n)+\sup\{\phi(\boldsymbol{\xi}):\boldsymbol{\xi}\notin\Omega_n\};$$ letting 
$t\rightarrow\infty$, yields $\varphi(\boldsymbol{\theta})\leq\sup_{\boldsymbol{\xi}\notin\Omega_n}\phi(\boldsymbol{\xi})\longrightarrow 0$, which proves (\ref{phi:eq}), hence (\ref{regio:eq}).

Finally, to prove continuity of $\phi$, note that both expressions in the r.h.s.~of (\ref{phi:eq}) define continuous functions, hence we only need to verify that $\varphi$ is continuous at boundary points of $\Theta_{\rm s}$. To this end, let $\boldsymbol{\theta}_*\in\partial\Theta_{\rm s}$. Since $\Theta_{\rm s}$ is open, we have $\partial\Theta_{\rm s}\subseteq\Theta_{\rm s}^\complement$, hence $\varphi(\boldsymbol{\theta}_*)=0$. On the other hand, $\varphi=\inf_t\varphi_t$ is upper semi-continuous, whence $$0\leq\underset{\boldsymbol{\theta}\rightarrow\boldsymbol{\theta}_*}{\lim\sup}\:
\varphi(\boldsymbol{\theta})\leq\varphi(\boldsymbol{\theta}_*)=0,$$
i.e., $\lim_{\boldsymbol{\theta}\rightarrow\boldsymbol{\theta}_*}\varphi(\boldsymbol{\theta})=0$. Therefore, $\varphi$ is continuous at $\boldsymbol{\theta}_*$, which proves the claim.

III. Let now $\mathbf{0}\leq\boldsymbol{\theta}<\boldsymbol{\vartheta}$ be s.t.\ $\boldsymbol{\vartheta}\in\Theta_{\rm s}$; in particular, we have $\boldsymbol{\theta}\in\Theta_{\rm s}$ (cf.\ II) and both $\pi_{\boldsymbol{\theta}}$ and $\pi_{\boldsymbol{\vartheta}}$ are supported on $\mathbb{X}_0$. 

First, one infers from (\ref{diff-int:eq}), (\ref{eq:aux}) and (\ref{embed:prob}) that
$$\forall t\geq 0:\:
\left[\exp(t{\boldsymbol A}_{\boldsymbol{\vartheta}})-\exp(t{\boldsymbol A}_{\boldsymbol{\theta}})\right]
\boldsymbol {Q}_{\boldsymbol{\theta}}\phi\leq\mathbf{0};$$
letting $t\rightarrow\infty$, yields $\langle\pi_{\boldsymbol{\vartheta}}-\pi_{\boldsymbol{\theta}},\boldsymbol {Q}_{\boldsymbol{\theta}}\phi\rangle\leq 0$.
Furthermore, using the identity $\pi_{\boldsymbol{\theta}}\boldsymbol{A}_{\boldsymbol{\theta}}=\mathbf{0}$ (valid for all $\boldsymbol{\theta}$'s) one obtains
\begin{eqnarray*}
\langle\pi_{\boldsymbol{\vartheta}}-\pi_{\boldsymbol{\theta}},\phi\rangle & = &
\langle\pi_{\boldsymbol{\vartheta}}-\pi_{\boldsymbol{\theta}},(\boldsymbol{Q}_{\boldsymbol{\theta}}-(1/a)
\boldsymbol {A}_{\boldsymbol{\theta}})\phi\rangle \\
& \leq & -(1/a)\cdot\langle\pi_{\boldsymbol{\vartheta}}-\pi_{\boldsymbol{\theta}},\boldsymbol {A}_{\boldsymbol{\theta}} \phi\rangle=
(1/a)\cdot\langle\pi_{\boldsymbol{\vartheta}},(\boldsymbol{A}_{\boldsymbol{\vartheta}}-\boldsymbol{A}_{\boldsymbol{\theta}})\phi\rangle \\
& = & (1/a)\sum_{k=1}^d(\vartheta_k-\theta_k)\cdot\langle\pi_{\boldsymbol{\vartheta}},(\Phi[f_{(0,k)}]-\boldsymbol {I}) \phi\rangle.
\end{eqnarray*}
Finally, we need to prove that the last expression in the above display is strictly negative for $\boldsymbol{\theta}<\boldsymbol{\vartheta}$.
To this end, let $k$ be such that $\theta_k<\vartheta_k$. Since $\pi_{\boldsymbol{\vartheta}}$ is supported on $\mathbb{X}_0$ it is enough to show that $(\Phi[f_{(0,k)}]-{\boldsymbol I})\phi$ may not vanish everywhere on $\mathbb{X}_0$. Indeed, let $\boldsymbol{\xi}_0:=\emptyset$ and $\boldsymbol{\xi}_{n+1}:=f_{(0,k)}(\boldsymbol{\xi}_n)$, for $n\geq 0$; note that $\boldsymbol{\xi}_n\in\mathbb{X}_0$, for all $n\geq 0$. Assuming that $(\Phi[f_{(0,k)}]-{\boldsymbol I})\phi$ vanishes on $\mathbb{X}_0$, it follows (by induction) that $\phi(\boldsymbol{\xi}_n)=\phi(\emptyset)>0$ (by assumption), for all $n\geq 0$, hence $\phi$ is constant and strictly positive along the infinite sequence $\{\boldsymbol{\xi}_n\}_{n\geq 0}\subseteq\mathbb{X}$. But since $\phi$ must vanish at infinity, this is a contradiction. This completes the proof of the theorem.
\hfill$\square$ 

\newpage

\section*{Acknowledgments}
The authors are grateful to A.M.\ Oprescu, researcher at University of Amsterdam, for assisting with the implementation of the numerical
experiments which led to establishing the results in this paper.

\begin{figure}[ht!]
\unitlength=4.8 mm 
\begin{center}
\begin{picture}(28.0,14.0)
\drawline(-0.40,0.0)(28.4,0.0)
\drawline(28.4,0.0)(28.4,14.0)
\drawline(28.4,14.0)(-0.40,14.0)
\drawline(-0.40,14.0)(-0.40,0.0)
\put(1.9,4.0){\makebox(0,0)[cc]{$\{\boldsymbol{\xi}_1''=\boldsymbol{\xi}_0''\}$}}
\put(1.9,10.0){\makebox(0,0)[cc]{$\{\boldsymbol{\xi}_1''\neq\boldsymbol{\xi}_0''\}$}}
\put(5.0,3.0){\makebox(0,0)[cc]{$\boldsymbol{\xi}$}}
\put(5.0,5.0){\makebox(0,0)[cc]{$\boldsymbol{\zeta}$}}
\put(5.6,3.0){\vector(1,3){2.00}}
\put(5.6,5.0){\vector(1,3){2.00}}
\textcolor{red}{\put(5.6,5.0){\vector(1,0){2.00}}}
\put(9.2,5.0){\makebox(0,0)[cc]{$\boldsymbol{\zeta}$}}
\put(9.2,9.0){\makebox(0,0)[cc]{$f_{(k,l)}(\boldsymbol{\xi})$}}
\put(9.2,11.0){\makebox(0,0)[cc]{$f_{(k,l)}(\boldsymbol{\zeta})$}}
{\put(10.8,11.0){\color{blue}\vector(1,0){14.4}}}
{\put(10.8,9.0)
{\color{blue}\vector(1,0){14.4}}}
\put(18.0,10.0){\makebox(0,0)[cc]{$(\mathfrak{H}_n)$}}
\put(26.4,11.0){\makebox(0,0)[cc]{$\boldsymbol{\xi}_{n+1}''$}}
\put(26.4,9.0){\makebox(0,0)[cc]{$\boldsymbol{\xi}_{n+1}'$}}
{\put(10.8,5.0){\color{blue}\vector(1,0){14.4}}}
\put(26.4,5.0){\makebox(0,0)[cc]{$\boldsymbol{\xi}_{n+1}''$}}
{\put(5.6,3.0){\color{blue}\vector(1,0){14.4}}}
\put(22.2,3.0){\makebox(0,0)[cc]{$\boldsymbol{\xi}_n'$}}
\textcolor{red}{\put(23.2,3.0){\vector(1,0){2.00}}}
\put(15.5,4.0){\makebox(0,0)[cc]{$(\mathfrak{H}_n)$}}
\dottedline(5.6,3.0)(10.8,5.0)
\dottedline(20.0,3.0)(25.2,5.0)
\put(26.4,3.0){\makebox(0,0)[cc]{$\boldsymbol{\xi}_{n+1}'$}}
\put(4.0,2.0){\vector(1,0){23.6}}
\drawline(4.0,2.0)(4.0,12.0)
\put(5.0,1.2){\makebox(0,0)[cc]{$0$}}
\put(9.2,1.2){\makebox(0,0)[cc]{$1$}}
\put(15.4,1.2){\makebox(0,0)[cc]{$\ldots$}}
\put(22.2,1.2){\makebox(0,0)[cc]{$n$}}
\put(26.4,1.2){\makebox(0,0)[cc]{$n+1$}}
\put(5.0,2.0){\line(0,-1){0.2}}
\put(9.2,2.0){\line(0,-1){0.2}}
\put(22.2,2.0){\line(0,-1){0.2}}
\put(26.4,2.0){\line(0,-1){0.2}}
\end{picture}\end{center}
\caption{Graphic representation of the proof (induction step) of Theorem \ref{HQ:them}. Pairs of parallel blue arrows represent the coupled sample paths presumed by the induction hypothesis $(\mathfrak{H}_n)$, whereas the red arrows represent one-step transitions conditioned to non-arrival events.}\label{proof:fig}
\end{figure}
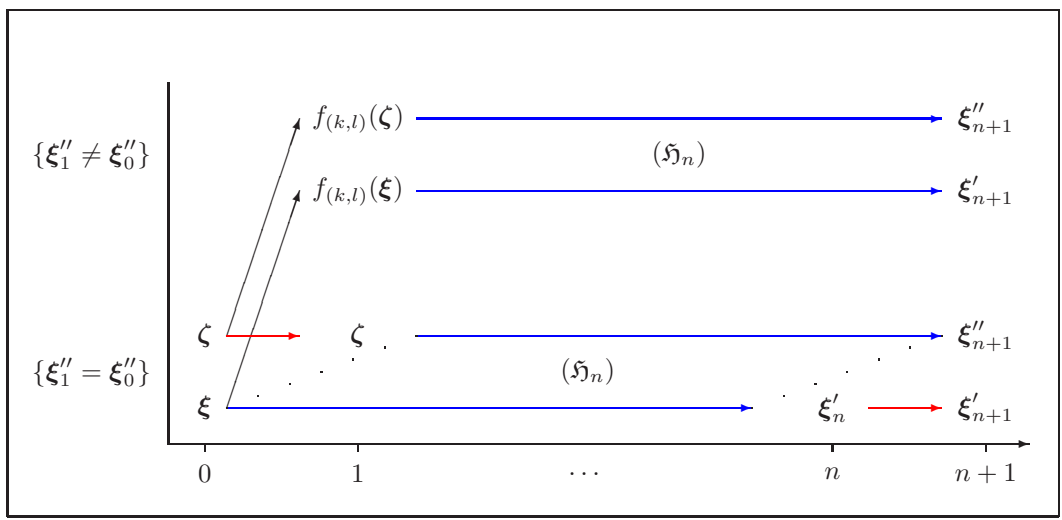

\newpage

\bibliographystyle{amsplain}
\bibliography{Biblist}

\providecommand{\bysame}{\leavevmode\hbox to3em{\hrulefill}\thinspace}
\providecommand{\MR}{\relax\ifhmode\unskip\space\fi MR }
\providecommand{\MRhref}[2]{%
  \href{http://www.ams.org/mathscinet-getitem?mr=#1}{#2}
}
\providecommand{\href}[2]{#2}
\begin{thebibliography}{10}

\bibitem{RS:93}
{A.N. Rybko} and {A.L. Stolyar}, \emph{Ergodicity of stochastic processes
  describing the operation of open queueing networks}, Problemy Peredachi
  Informatsii \textbf{28} (1993), pp. 3--26.

\bibitem{Kelly}
{F.P. Kelly}, \emph{Networks of queues with customers of different types},
  Journal of Applied Probability \textbf{12} (1975), pp. 542--554.

\bibitem{LM:2016}
{H. Leahu} and {M. Mandjes}, \emph{A numerical approach to stability of
  multiclass queueing networks}, IEEE Transactions on Automatic Control
  \textbf{62} (2017), pp. 5478--5484.

\bibitem{Dai:95}
{J.G. Dai}, \emph{On positive {H}arris recurrence of multiclass queueing
  networks: A unified approach via fluid limit models}, Annals of Applied
  Probability \textbf{5} (1995), pp. 49--77.

\bibitem{Dai:99}
{J.G. Dai}, {J.J. Hasenbein}, and {J.H. Vande Vate}, \emph{Stability of a
  three-station fluid network}, Queueing Systems \textbf{33} (1999), 293--325.

\bibitem{Dai:04}
\bysame, \emph{Stability and instability of a two-station queueing network},
  Annals of Applied Probability \textbf{14} (2004), 326--377.

\bibitem{S-Yao}
{J.G. Shanthikumar} and {D.D. Yao}, \emph{Stochastic monotonicity in general
  queueing networks}, Journal of Applied Probability \textbf{26} (1989), pp.
  413--417.

\bibitem{Bramson:94}
{M. Bramson}, \emph{Instability of {F}{I}{F}{O} queueing networks}, Annals of
  Applied Probability \textbf{4} (1994), no.~2, pp. 414--431.

\bibitem{Bramson:95}
\bysame, \emph{Instability of {F}{I}{F}{O} queueing networks with quick service
  times}, Annals of Applied Probability (1994), no.~3, 693--718.

\bibitem{Bramson:08}
\bysame, \emph{Stability of queueing networks}, Probability Surveys \textbf{5}
  (2008), pp. 169--345.

\bibitem{KS:90}
{P.R. Kumar} and {T.I. Seidman}, \emph{Dynamic instabilities and stabilization
  methods in distributed real-time scheduling of manufacturing systems}, IEEE
  Transactions on Automatic Control \textbf{35} (1990), pp. 289--298.

\bibitem{AG:2007}
{S. Asmussen} and {P.W. Glynn}, \emph{Stochastic simulation: Algorithms and
  analysis}, Springer-Verlag NY, 2007.

\bibitem{MT:93}
{S.P. Meyn} and {R.L. Tweedie}, \emph{Markov chains and stochastic stability},
  Springer Verlag, London, 1993.

\bibitem{Meyn-Tweedie:II}
\bysame, \emph{Stability of {M}arkovian processes {II}: Continuous-time
  processes and sampled chains}, Advances in Applied Probability \textbf{25}
  (1993), pp. 487--517.

\bibitem{S:94}
{T.I. Seidman}, \emph{`{F}irst come, first serve' can be unstable!}, IEEE
  Transactions on Automatic Control \textbf{39} (1994), pp. 2166--2171.

\bibitem{Dumas}
{V. Dumas}, \emph{A multiclass network with non-linear, non-convex,
  non-monotonic stability conditions}, Queueing Systems \textbf{25} (1997),
  1--43.

\bibitem{Whitt:93}
{W. Whitt}, \emph{Large fluctuations in a deterministic multiclass network of
  queues}, Management Science \textbf{39} (1993), 1020--1028.

\bibitem{Massey:87}
{W.A. Massey}, \emph{Stochastic orderings for {M}arkov processes on partially
  ordered spaces}, Mathematics of Operations Research \textbf{12} (1987), pp.
  350--367.

\end{thebibliography}
\end{document}